\documentclass[11pt]{amsart}
\usepackage{amsmath, amssymb, amsthm}
\usepackage{amsrefs}
\usepackage{graphicx}
\usepackage{stmaryrd}
\usepackage{url}
\usepackage[most]{tcolorbox}
\usepackage{xcolor}
\usepackage{caption}

\makeatletter
\renewcommand{\section}{\@startsection{section}{1}\z@{2.2\linespacing\@plus\linespacing}{1.1\linespacing}{\normalfont\scshape\centering}}
\renewcommand{\subsection}{\@startsection{subsection}{2}\z@{1.1\linespacing\@plus.7\linespacing}{-.5em}{\normalfont\bfseries}}
\makeatother

\def\ucr#1{{\text{\rm cr}(#1)}}
\def\iocr#1{{\text{\rm iocr}(#1)}}
\def\ocr#1{{\text{\rm ocr}_{+}(#1)}}
\def\pcr#1{{\text{\rm pcr}_{+}(#1)}}
\def\floor#1{{\lfloor{#1}\rfloor}}
\def\Bigfloor#1{{\Big\lfloor{#1}\Big\rfloor}}
\def\real{{\mathbb R}}
\def\ss{{\mathcal S}}
\def\xx{{\mathcal X}}
\def\pp{{\mathcal P}}
\def\Cr{{\text{\rm Quad}}}
\def\Rot{{\text{\rm Rot}}}
\def\sign{{\text{\rm sign}}}
\def\nr{{\textnormal{\textsc{nr}}}}
\def\crq#1{{\boxtimes(#1)}}
\def\orient#1{{\circlearrowleft_{#1}}}

\tcbset{before skip=18pt,after skip=18pt}
\newtcolorbox{defbox}[1]{colback=gray!4,colframe=gray!50,colbacktitle=gray!22,coltitle=black!85,boxrule=0.6pt,arc=2pt,left=8pt,right=8pt,top=5pt,bottom=5pt,fonttitle=\bfseries,title={#1}}
\newtcolorbox{settingsbox}[1]{colback=gray!4,colframe=gray!50,colbacktitle=gray!22,coltitle=black!85,boxrule=0.6pt,arc=2pt,left=8pt,right=8pt,top=12pt,bottom=12pt,fonttitle=\bfseries,title={#1}}

\newtheorem{theorem}{Theorem}
\newtheorem{lemma}[theorem]{Lemma}
\newtheorem{corollary}[theorem]{Corollary}
\newtheorem{proposition}[theorem]{Proposition}
\newtheorem{conjecture}[theorem]{Conjecture}
\newtheorem{remark}[theorem]{Remark}
\newtheorem{observation}[theorem]{Observation}

\theoremstyle{definition}
\newtheorem{definition}[theorem]{Definition}

\begin{document}


\author{Gelasio Salazar}
\address{Instituto de F\'\i sica, Universidad Aut\'onoma de San Luis Potos\'{\i}. San Luis Potos\'{\i}, Mexico.}
\email{\tt gelasio.salazar@uaslp.mx}

\title[The Harary--Hill conjecture for rotation systems]{The proof of the Harary--Hill conjecture: a detailed exposition and an extension to rotation systems}

\begin{abstract}
The Harary--Hill conjecture states that the crossing number $\ucr{K_n}$ of the complete graph $K_n$ is $Z(n):=\frac{1}{4}\lfloor\frac{n}{2}\rfloor\lfloor\frac{n-1}{2}\rfloor\lfloor\frac{n-2}{2}\rfloor\lfloor\frac{n-3}{2}\rfloor$ for all $n\ge 3$. In October 2026 OpenAI uploaded to GitHub a proof of this conjecture. The proof is very short, but it relies on polynomial methods that are not part of the standard techniques used by many combinatorialists, and it is written with very little motivation or explanation. Since we found the proof very hard to follow, we decided to write this paper, which contains a complete proof of the conjecture written for mathematicians who are more familiar with combinatorics than with polynomial methods. In fact we prove a more general result on rotation systems of $K_n$, from which the Harary--Hill conjecture easily follows. Namely, we show that for every rotation system $\Pi$ of $K_n$ the number of $4$-vertex subsystems of $\Pi$ realizable by drawings of $K_4$ with one crossing, plus three times the number of non-realizable $4$-vertex subsystems of $\Pi$, is at least $Z(n)$. The statement and proof of our result are purely combinatorial and do not involve drawings at all.
\end{abstract}

\date{\today}

\subjclass[2020]{05C10, 05C62}
\keywords{crossing number, complete graph, Harary--Hill conjecture, rotation system, polynomial method}

\maketitle


\section{Introduction}\label{sec:intro}

In a {\em drawing} of a graph in the plane, vertices are represented by distinct points and each edge is represented by a Jordan arc whose endpoints are the endvertices of the edge. No edge passes through a vertex, and edges intersect each other in a finite number of points, each of which is either a common endvertex or a point where the two edges cross each other. A drawing is {\em simple} if no two edges share more than one point (either a common endvertex or a proper {\em crossing}), and no three edges meet at a common interior point. 

The {\em crossing number} $\ucr{D}$ of a drawing $D$ is the number of pairwise crossings of edges in $D$. The {\em crossing number} $\ucr{G}$ of a graph $G$ is the minimum value of $\ucr{D}$ taken over all drawings of $G$ in the plane. Every graph $G$ has a simple drawing $D$ such that $\ucr{D}=\ucr{G}$~\cite{Kainen68,schaefer,Szekely04}, and so for the purposes of computing crossing numbers it suffices to consider simple drawings. For this reason throughout this paper all drawings under consideration are assumed to be simple.


\subsection{The Harary--Hill conjecture and the dual purpose of this paper}\label{sub:purpose}

For a long time one of the most important open problems in topological graph theory was the Harary--Hill conjecture~\cite{HararyHill}, which states the following.\footnote{This is also known as {\em Guy's conjecture}, since it first appeared in a paper by Richard K.~Guy in 1960~\cite{Guy60}.}

\begin{conjecture}[Harary--Hill, 1963]\label{con:hill}
For every integer $n\ge 3$
$$\ucr{K_n} = Z(n):=\frac{1}{4}\Bigfloor{\frac{n}{2}}\Bigfloor{\frac{n-1}{2}}\Bigfloor{\frac{n-2}{2}}\Bigfloor{\frac{n-3}{2}}.$$
\end{conjecture}

A well-known straightforward counting argument shows the following (see for instance~\cite{szekely}).

\begin{observation}\label{obs:oddn}
If the Harary--Hill conjecture holds for some odd integer $n\ge 3$, then it also holds for $n+1$.
\end{observation}

In October 2026 OpenAI uploaded to GitHub a proof of the following theorem. In view of Observation~\ref{obs:oddn}, this proves Conjecture~\ref{con:hill}.

\begin{theorem}[OpenAI~\cite{openai}. Implies Conjecture~\ref{con:hill}]\label{thm:hill}
For every odd integer $n\ge 3$, $\ucr{K_n}=Z(n)$.
\end{theorem}

We have been deeply interested in this conjecture for a long time, and we were highly motivated to read the proof. Since we found it very hard to follow, we decided to write an expository paper in order to understand the main ideas of the proof. 

After writing a first draft we realized that the tools developed in~\cite{openai} could be adjusted to prove a strictly more general result, and so we moved on to write this paper, which maintains the original spirit of a detailed exposition but also contains a complete proof of this stronger statement. We emphasize that this paper is written having as target audience a combinatorialist (such as myself) with little or no knowledge about the polynomial methods that are at the heart of the proof. 


\subsection{Our main result}\label{sub:ourmain}

We actually extend the techniques in~\cite{openai} to prove a stronger statement than Theorem~\ref{thm:hill}, in terms of rotation systems rather than drawings. The notions of a rotation system and of a realizable and non-realizable rotation system are reviewed in Section~\ref{sub:rotsys}. In the meantime we state our main result, for the benefit of the reader familiar with rotation systems.

In the statement we use the following two numbers associated to (abstract) rotation systems of $K_n$. We recall that all drawings under consideration are simple. Given a rotation system $\Pi$ of $K_n$:

\begin{itemize}

\item we let $\crq{\Pi}$ denote the number of rotation subsystems of $\Pi$ of size $4$ that are realizable by drawings of $K_4$ with one crossing, and

\item we let $\nr(\Pi)$ denote the number of rotation subsystems of $\Pi$ of size $4$ that are non-realizable by drawings of $K_4$.

\end{itemize}

\begin{theorem}[Implies Conjecture~\ref{con:hill}]\label{thm:odd}
Let $\Pi$ be any rotation system of $K_n$, where $n\ge 3$ is an integer. Then
$$
\crq{\Pi}+3\,\nr(\Pi)\ge Z(n).
$$
\end{theorem}

As we discuss in Section~\ref{sec:concluding} at the end of the paper, it seems plausible that this inequality holds even if the coefficient of $\nr(\Pi)$ is $1$ rather than $3$, and this would establish an appealing conjecture formulated by Dan Archdeacon in the 1990s.

A reader familiar with the Harary--Hill conjecture and with rotation systems will have no trouble seeing that Theorem~\ref{thm:odd} easily implies the Harary--Hill conjecture. Indeed, let $n\ge 3$ be an integer, let $D$ be a drawing of $K_n$, and let $\Pi$ be its rotation system. Since $\Pi$ arises from $D$, none of the rotation subsystems of $\Pi$ of size $4$ is non-realizable, and so $\nr(\Pi)=0$. Therefore in this case Theorem~\ref{thm:odd} says that $\crq{\Pi}\ge Z(n)$. Since each crossing of $D$ lies in exactly one of the drawings of $K_4$ induced by $D$, the number $\crq{\Pi}$ is the number of crossings of $D$. Thus $D$ has at least $Z(n)$ crossings, and so $\ucr{K_n}\ge Z(n)$ (the difficult part of the equality in Conjecture~\ref{con:hill}) follows.

We conclude this introduction with some background on the conjecture and an exhaustive discussion on previous work on this problem.


\subsection{Background on the Harary--Hill conjecture}\label{sub:background}

As Beineke and Wilson explain in their lively survey~\cite{beineke-wilson}, the problem of estimating the crossing number of complete graphs seems to have been first explored by the British artist Anthony Hill in the late 1950s. Hill found a way to draw $K_n$ with $Z(n)$ crossings for every integer $n \ge 3$, and later approached graph theorist Frank Harary, and together they wrote a paper where they stated Conjecture~\ref{con:hill} above~\cite{HararyHill}. 

The conjectured value was first published by Guy~\cite{Guy60}. Saaty stated the same conjecture independently in 1964, and argued that the conjecture would follow if some crossing-minimal drawing of $K_n$ contained crossing-minimal drawings of $K_{n-2}$, $K_{n-4}$, and $K_{n-6}$~\cite{Saaty,Saaty67}. The problem was popularized by Erd\H{o}s and Guy in their classical paper on crossing number problems~\cite{ErdosGuy}. We also refer the reader to the survey by Sz\'ekely~\cite{szekely}, and to Schaefer's dynamic survey on crossing numbers~\cite{SchaeferSurvey}.

\subsubsection{Constructions of drawings of $K_n$ with $Z(n)$ crossings}
\begin{figure}[ht]
\centering
\includegraphics{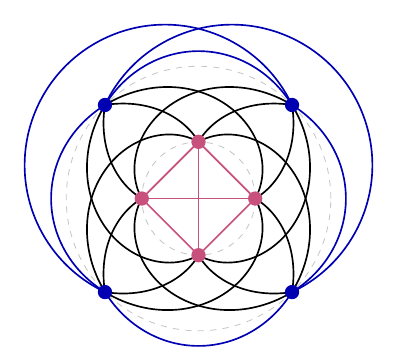}
\caption{Hill's drawing of $K_8$. Four vertices are drawn on an inner circle (dashed) and four vertices are drawn on an outer circle (also dashed). Edges joining two inner vertices stay inside the inner circle, edges joining two outer vertices stay outside the outer circle, and each inner vertex is joined to each outer vertex with an edge drawn inside the annulus bounded by the inner and outer circles. The drawing has $Z(8)=18$ crossings.}
\label{fig:tincan}
\end{figure}

There are many distinct constructions that yield drawings of $K_n$ with exactly $Z(n)$ crossings. Hill's construction produces {\em cylindrical} (also known as {\em tin can}) drawings, which are drawings in which the vertices are drawn on two concentric circles, and no edge intersects any of these circles, except at its endpoints. We refer the reader to Figure~\ref{fig:tincan} for an illustration of such a cylindrical drawing of $K_8$.

Bla\v zek and Koman~\cite{BlazekKoman} found an alternative construction with $Z(n)$ crossings of $K_n$. Their construction yields $2$-{\em page drawings} of $K_n$, that is, drawings in which every vertex lies on the $x$-axis, and each edge lies (except for its endpoints) either in the upper or in the lower halfplane. \'{A}brego, Aichholzer, Fern\'{a}ndez-Merchant, Ramos, and Vogtenhuber~\cite{NonShellable} described a variant of Hill's construction that yields drawings of $K_n$ with $Z(n)$ crossings, for every odd $n\ge 11$. We refer the reader to Figures 13.3 and 13.4 in~\cite{szekely} for lively descriptions of the cylindrical drawings devised by Hill and the $2$-page drawings of Bla\v zek and Koman.

Further constructions have been found recently. Mohar~\cite{Mohar} gave a very general construction of drawings of $K_n$ with $Z(n)$ crossings, which also gives a simple explanation of Moon's early result~\cite{Moon} that the expected number of crossings in a random spherical geodesic drawing of $K_n$ is $(1/4)(n/2)((n-1)/2)((n-2)/2)((n-3)/2)$, eerily close to $Z(n)$. Felsner~\cite{Felsner} recently reported new constructions of drawings of $K_n$ with $Z(n)$ crossings.

In a MathOverflow entry dating back to 2013 Kyn\v{c}l~\cite{kynclmathov} asked about drawings of $K_n$ with $Z(n)$ crossings, different from the constructions that were known back then. In a very recent answer to this question~\cite{kynclmathov2026} he pointed out that there are at least $n^{2n+o(n)}$ labeled spherical geodesic drawings of $K_n$ with $Z(n)$ crossings. Kyn\v{c}l derived this bound noting that Wagner~\cite[Proposition 1.4]{wagner2006},~\cite[Observation 2.6, Proposition 4.1, Observation 9.1 and Corollary 9.2]{wagner2008} observed that by Gale's duality, spherical drawings of $K_n$ with $Z(n)$ crossings correspond to neighborly polytopes with $n$ vertices in $\real^{n-4}$, and there are constructions of $n^{2n+o(n)}$ such labeled neighborly polytopes~\cite{padrol,gonskapadrol}.

In another edit from 2013 to the same MathOverflow entry Kyn\v{c}l introduced the notion of {\em pseudospherical drawings}, which generalize spherical geodesic drawings in a similar way as pseudolinear drawings generalize rectilinear drawings. Pseudospherical drawings were further investigated more recently by Arroyo, Richter, and Sunohara~\cite{ARS}.

In the same recent answer mentioned above Kyn\v{c}l~\cite{kynclmathov2026} also notes that there are $2^{\Omega(n^2)}$ pseudospherical drawings of $K_n$ with $Z(n)$ crossings. We refer the reader to~\cite{kynclmathov} for details on how this bound is obtained, and only mention that it follows by a construction by Kortenkamp~\cite{kortenkamp} on spherical drawings of $K_{2n}$ with $Z(2n)$ crossings, which actually holds in the more general setting of oriented matroids, which in turn correspond to pseudospherical drawings.

\subsubsection{Families of drawings for which the Harary--Hill conjecture had been verified}

The Harary--Hill conjecture was verified by \'Abrego, Aichholzer, Fern\'andez-Merchant, Ramos, and Salazar both for $2$-page~\cite{2Page} and for cylindrical~\cite{Cylindrical} drawings. It is also known that the conjecture holds for {\em monotone} drawings, that is, drawings in which each edge is drawn as an $x$-monotone curve~\cite{Monotone2,BalkoMonotone}. The alternative construction in~\cite{NonShellable} yields drawings that are neither $2$-page nor cylindrical, but they satisfy a property called {\em bishellability}. In~\cite{Bishellable}, it was proved that the Harary--Hill conjecture holds for bishellable drawings (more precisely, for $(\lfloor n/2\rfloor-2)$-bishellable drawings, see~\cite{Bishellable} for the definition).

Mutzel and Oettershagen~\cite{MutzelOettershagen} went one step further, introducing the class of seq-shellable drawings, which strictly contains the class of bishellable drawings, and proving the Harary--Hill conjecture for this larger class. In a related work Arroyo, McQuillan, Richter, and Salazar~\cite{Convex} introduced convex drawings of $K_n$, a class that lies between rectilinear and general drawings, and discussed properties of optimal drawings of $K_n$.

In its full generality (prior to~\cite{openai}), the conjecture had only been verified for $n\le 14$. For $n\le 10$ this was first reported by Saaty~\cite{Saaty69} in 1969 (see also~\cite{Guy10,Guy72}). Using these values, Saaty also showed~\cite{Saaty71} that if $\ucr{K_n}$ is given by one polynomial in $n$ for even $n$ and by another one for odd $n$, then these two polynomials are the ones that give $Z(n)$. More recently, McQuillan and Richter~\cite{DanBruce} gave a computer-free verification of the Harary--Hill conjecture for $n=9$ (from which it follows for $n=10$, by Observation~\ref{obs:oddn}). Pan and Richter~\cite{PanRichter} gave a computer-assisted proof for $n=11$ (and hence for $n=12$). The Harary--Hill conjecture for $n\le 12$ has also been verified in~\cite{EuroCG15}.

McQuillan, Pan, and Richter~\cite{K13} ruled out the possibility that $\ucr{K_{13}}=217$, and since $\ucr{K_{13}}$ is an odd number~\cite{Kleitman,Parity}, it followed that $\ucr{K_{13}}\in \{219,221,223,225\}$. In~\cite{EuroCG15} this was further narrowed to $\ucr{K_{13}} \in \{223,225\}$. Finally, in 2021 Aichholzer~\cite{Aichholzer13} announced the proof of $\ucr{K_{13}}=Z(13)=225$ (and hence of $\ucr{K_{14}}=Z(14)$).

We note that in a recent major development Streltsova and Wagner~\cite{StreltsovaWagner} proved that the Harary--Hill conjecture holds for spherical geodesic drawings. More precisely, if $n$ points are placed in general position on the sphere and every pair is joined by its unique shortest geodesic arc, then the resulting drawing has at least $Z(n)$ crossings.

Two very recent results also involve verifications of the conjecture. Felsner~\cite{Felsner} proved the conjecture for antipodally shellable drawings, and Hubard, Ndiaye, and Sol\'e-Pi~\cite{HubardNdiayeSolePi} gave a short proof for antipodal drawings, in a paper that also proposes a generalization of the conjecture to simplicial complexes.

\subsubsection{Relationship with Zarankiewicz's conjecture}

Besides the Harary--Hill conjecture the most famous problem involving crossing numbers is Zarankiewicz's conjecture, which claims that the crossing number of the complete bipartite graph $K_{m,n}$ is $\ucr{K_{m,n}}=Z(m,n):=\floor{m/2}\floor{(m-1)/2}\floor{n/2}\floor{(n-1)/2}$ for all $m,n\ge 3$. Zarankiewicz~\cite{Zarankiewicz} gave drawings of $K_{m,n}$ with $Z(m,n)$ crossings and claimed a proof that they are optimal (hence the name of the conjecture). However, Kainen and Ringel independently found out that his proof had an unpluggable gap~\citelist{\cite{Guy69}}. This famous problem goes back to Tur\'an, who explains that the question had its origin in a brick factory during the Second World War~\cite{TuranWelcome} (hence the title of the survey by Beineke and Wilson~\cite{beineke-wilson}). 

Kainen~\cite{Kainen68} proved that Zarankiewicz's conjecture implies an asymptotic version of the Harary--Hill conjecture, namely that $\ucr{K_n}/Z(n)\to 1$ as $n\to\infty$. Richter and Thomassen~\cite{RichterThomassen} showed that for this conclusion it suffices to assume the asymptotic version of Zarankiewicz's conjecture for the balanced graphs $K_{n,n}$. More precisely, they proved that the limit of $\ucr{K_n}/Z(n)$ exists and is at least the limit of $\ucr{K_{n,n}}/Z(n,n)$. 

We note that at the same time that OpenAI announced the proof of the Harary--Hill conjecture it also announced a proof of Zarankiewicz's conjecture~\cite{openai-bip}.

\subsubsection{General lower bounds for $\ucr{K_n}$}

For a long time the best lower bound for $\ucr{K_n}$ was the asymptotic bound $\ucr{K_n} \ge 0.8Z(n)$, which follows from an elementary counting argument due to Kainen~\cite{Kainen68} combined with Kleitman's result~\cite{Kleitman70} that $\ucr{K_{5,n}}=4\floor{n/2}\floor{(n-1)/2}$, confirming that Zarankiewicz's conjecture holds for $m=5$. 

More recent lower bounds on $\ucr{K_n}$ were obtained by computer-assisted methods. De Klerk, Maharry, Pasechnik, Richter, and Salazar~\cite{dKMPRS} used semidefinite programming to obtain a lower bound for $\ucr{K_{7,n}}$, from which it follows that $\ucr{K_n}\ge 0.83\,Z(n)$ for all sufficiently large $n$. De Klerk, Pasechnik, and Schrijver~\cite{dKPS} later improved this to $\ucr{K_n}\ge 0.859\,Z(n)$. Using flag algebras, Norin and Zwols~\cite{NorinZwols} obtained an asymptotic lower bound for the crossing number of complete bipartite graphs, from which it follows that $\ucr{K_n}>0.905\,Z(n)$ for all sufficiently large $n$. 

Balogh, Lidick\'y, and Salazar~\cite{BLS} used flag algebras to prove that $\ucr{K_n}\ge 0.985\,Z(n)$ for all sufficiently large $n$. This was recently improved by Balogh, Layman, Lidick\'y, Pfender, Salazar, Spiro, and Volec~\cite{flagzar}, who showed that $\ucr{K_n} \ge 0.9899\, Z(n)$ for all sufficiently large $n$. Using Laci Sz\'ekely's lively expression, this last result can be paraphrased by saying that the Harary--Hill conjecture is asymptotically ``at least $98.99\%$ true''. 

We finally mention that the OpenAI proof uses algebraic information about a drawing, namely the crossings between two closed walks counted with signs that record the orientation of each crossing.  Such signed crossing counts of pairs of cycles go back to Tutte's algebraic theory of crossing numbers~\cite{Tutte}.



\section{Our main result}\label{sec:main}

As we briefly mentioned in Section~\ref{sub:ourmain}, using the techniques developed in~\cite{openai} to prove the Harary--Hill conjecture we can actually establish a stronger statement, in terms of rotation systems rather than drawings. Its statement and proof involve no drawings at all. 

 In Section~\ref{sub:rotsys} we recall the notion of a rotation system and its subsystems. In Section~\ref{sub:mainresult} we state this main result (namely Theorem~\ref{thm:odd}) again, and show that it easily implies the Harary--Hill conjecture.  

Before getting into rotation systems let us lay out a convention for the rest of this paper.

\begin{settingsbox}{The vertex set of $K_n$ is $[n]$}
For the rest of this paper the vertices of $K_n$ are labelled $1,\ldots,n$, and we let $[n]:=\{1,\ldots,n\}$.
\end{settingsbox}

\subsection{Rotation systems and subsystems}\label{sub:rotsys}

Let us start with the central notion of this paper, namely a rotation system. Even though (as we will see shortly) rotation systems are purely combinatorial objects that can be (and will also be) defined without any reference to drawings, they naturally arise from drawings.

\begin{defbox}{The rotation system of a drawing}
\begin{minipage}[c]{0.6\linewidth}
Let $D$ be a drawing of $K_n$. For each vertex $i$, the {\em rotation} $\pi(i)$ at $i$ is the cyclic permutation that records the anticlockwise cyclic order in which the edges incident with $i$ leave $i$ in $D$, each edge named by its endvertex other than $i$. The collection $\Pi=(\pi(i))_{i\in[n]}$ is the {\em rotation system} of $D$.
\end{minipage}\hfill
\begin{minipage}[c]{0.36\linewidth}
\centering
\includegraphics{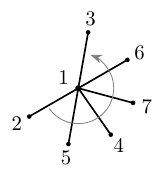}
\end{minipage}

\end{defbox}

\begin{defbox}{Convention: we write rotations as \underline{linear} permutations}
For definiteness we will write each rotation as a linear permutation, starting with its smallest element. Thus $\pi(1)=(2,\ldots)$, and $\pi(i)=(1,\ldots)$ for every $i=2,\ldots,n$. 
\medskip

We emphasize that a rotation is still a cyclic order. This convention, which turns out to be very convenient, only fixes where we start writing this cyclic order.
\end{defbox}

 For instance, in the illustration above, the rotation at $1$ is simply the linear permutation $\pi(1)=(2,5,4,7,6,3)$.

\begin{figure}[ht!]
\centering
\includegraphics{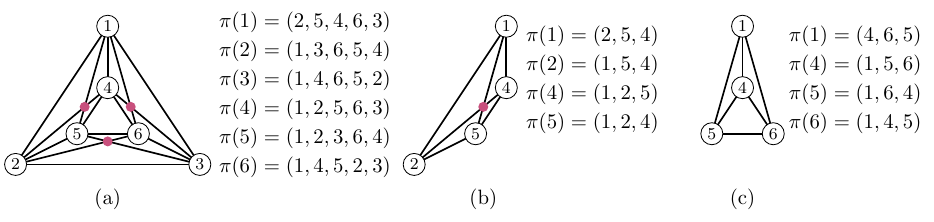}
\caption{In (a) we illustrate a drawing $D$ of $K_6$ with $Z(6)=3$ crossings and display its rotation system $\Pi=\bigl(\pi(1),\pi(2),\pi(3),\pi(4),\pi(5),\pi(6)\bigr)$. In (b) we show the one-crossing drawing of $K_4$ induced by $D$ on the vertices $1,2,4,5$, together with its rotation system. In (c) we have the plane drawing of $K_4$ induced by $D$ on the vertices $1,4,5,6$.}
\label{fig:k6}
\end{figure}

Figure~\ref{fig:k6}(a) shows a drawing of $K_6$ with three crossings, together with its rotation system.

The rotation system of a drawing of $K_n$ contains a great deal of information about the drawing. Gioan~\cite{Gioan05,Gioan22} proved that two drawings of $K_n$ in the sphere with the same rotation system are the same up to homeomorphism and a sequence of triangle flips. Other proofs of Gioan's theorem are given in~\cite{AMRS} and~\cite{SchaeferDetour}, and the latter is particularly short. A similar result for complete multipartite graphs, using a suitable extension of rotation systems, was proved by Aichholzer et al.~\cite{AichholzerEtAl23}. It is easy to see that the rotation system of a drawing of $K_n$ determines which pairs of edges cross, since this can be checked on each set of four vertices (see~\cite[Observation~14]{kyncl13}). Conversely (see for instance Gioan~\cite{Gioan05} and Kyn\v{c}l~\cite{kyncl}) the pairs of crossing edges determine the rotation system, up to reversing all the rotations.

As we mentioned above, a rotation system is a purely combinatorial entity, and so we may speak of a rotation system without any reference to a drawing. Let us formalize this.

\begin{defbox}{(Abstract) Rotation systems}
A {\em rotation system {of $K_n$}} is any collection $\Pi=(\pi(i))_{i\in[n]}$, where $\pi(i)$ is a cyclic permutation of $[n]\setminus\{i\}$, for each vertex $i$.

\medskip
More generally, if $S$ is a set of vertices of $K_n$, a {\em rotation system on $S$} assigns to each vertex of $S$ a cyclic permutation of the other vertices of $S$. That is, it is a rotation system of the complete graph with vertex set $S$.
\end{defbox}

Not every rotation system comes from a drawing. This leads to the following standard notion. We recall that all drawings under consideration are simple.

\begin{defbox}{Realizable and non-realizable rotation systems}
A rotation system of $K_n$ is {\em realizable} if it is the rotation system of a drawing of $K_n$. More generally, a rotation system on a set $S$ of vertices is {\em realizable} if it is the rotation system of a drawing of the complete graph with vertex set $S$. A rotation system that is not realizable is {\em non-realizable}.
\end{defbox}

The question of deciding whether a rotation system is realizable dates back at least to the 1990s~\cite{Archdeacon99}. Kyn\v{c}l~\cite{kyncl} gave a polynomial-time algorithm that decides whether a set of pairs of edges of $K_n$ is the set of pairs of crossing edges of a drawing of $K_n$. He later simplified it~\cite{kyncl20}, showing that the sets of pairs of crossing edges in drawings of $K_n$ are characterized by a finite list of forbidden configurations, each with at most six vertices. Combined with computations by \'Abrego et al.~\cite{EuroCG15}, this implies that for $n\ge 5$ a rotation system of $K_n$ is realizable if and only if all its subsystems of size $5$ are realizable~\cite{kyncl20}. Thus it follows that there is an $O(n^5)$ time algorithm that decides whether a rotation system is realizable.

\begin{defbox}{Rotation subsystems}
If $\Pi$ is a rotation system of $K_n$ and $S\subseteq [n]$, then $\Pi$ naturally induces a rotation system on $S$, simply by deleting from the rotation of each vertex of $S$ the vertices not in $S$. We use $\Pi|_{S}$ to denote this rotation system, and say that $\Pi|_{S}$ is a {\em subsystem} of $\Pi$ of size $|S|$.
\end{defbox}

We are particularly interested in the case in which $S$ is a set of $4$ vertices of $K_n$. These sets are so important that we coin a term for them.

\begin{defbox}{Quads}
A {\em quad} is a set of four vertices of $K_n$.
\end{defbox}

A rotation system on a quad $Q=\{i,j,k,\ell\}$, with $i<j<k<\ell$, consists of four rotations, and each of them is one of the two cyclic orders of the other three vertices of $Q$. Thus there are $2^4=16$ rotation systems on $Q$.

These $16$ rotation systems are listed in Table~\ref{tab:sixteen}, split into three groups, or {\em types}. The next observation explains where the three type labels come from.

\begin{table}[h!]
\centering
\small
\renewcommand{\arraystretch}{1.05}\setlength{\tabcolsep}{6pt}
\begin{tabular}{r|cccc|l}
 & $\pi(i)$ & $\pi(j)$ & $\pi(k)$ & $\pi(\ell)$ & type \\ \hline
1 & $(j,k,\ell)$ & $(i,\ell,k)$ & $(i,j,\ell)$ & $(i,k,j)$ & plane\\
2 & $(j,\ell,k)$ & $(i,k,\ell)$ & $(i,\ell,j)$ & $(i,j,k)$ & plane\\ \hline
3 & $(j,k,\ell)$ & $(i,\ell,k)$ & $(i,\ell,j)$ & $(i,j,k)$ & crossing\\
4 & $(j,\ell,k)$ & $(i,k,\ell)$ & $(i,j,\ell)$ & $(i,k,j)$ & crossing\\
5 & $(j,k,\ell)$ & $(i,k,\ell)$ & $(i,j,\ell)$ & $(i,j,k)$ & crossing\\
6 & $(j,\ell,k)$ & $(i,\ell,k)$ & $(i,\ell,j)$ & $(i,k,j)$ & crossing\\
7 & $(j,k,\ell)$ & $(i,k,\ell)$ & $(i,\ell,j)$ & $(i,k,j)$ & crossing\\
8 & $(j,\ell,k)$ & $(i,\ell,k)$ & $(i,j,\ell)$ & $(i,j,k)$ & crossing\\ \hline
9 & $(j,k,\ell)$ & $(i,k,\ell)$ & $(i,j,\ell)$ & $(i,k,j)$ & non-realizable\\
10 & $(j,k,\ell)$ & $(i,k,\ell)$ & $(i,\ell,j)$ & $(i,j,k)$ & non-realizable\\
11 & $(j,k,\ell)$ & $(i,\ell,k)$ & $(i,j,\ell)$ & $(i,j,k)$ & non-realizable\\
12 & $(j,k,\ell)$ & $(i,\ell,k)$ & $(i,\ell,j)$ & $(i,k,j)$ & non-realizable\\
13 & $(j,\ell,k)$ & $(i,k,\ell)$ & $(i,j,\ell)$ & $(i,j,k)$ & non-realizable\\
14 & $(j,\ell,k)$ & $(i,k,\ell)$ & $(i,\ell,j)$ & $(i,k,j)$ & non-realizable\\
15 & $(j,\ell,k)$ & $(i,\ell,k)$ & $(i,j,\ell)$ & $(i,k,j)$ & non-realizable\\
16 & $(j,\ell,k)$ & $(i,\ell,k)$ & $(i,\ell,j)$ & $(i,j,k)$ & non-realizable\\
\end{tabular}
\caption{The sixteen rotation systems on a quad $Q=\{i,j,k,\ell\}$ with $i<j<k<\ell$. By convention, each rotation is written starting at the smallest of the other three vertices of $Q$.}
\label{tab:sixteen}
\end{table}

\begin{observation}[Classification of rotation systems of size $4$]\label{obs:k4rows}
Let $\Pi$ be a rotation system on a quad $Q=\{i,j,k,\ell\}$, with $i<j<k<\ell$. Then $\Pi$ is one of the $16$ systems in Table~\ref{tab:sixteen}. Moreover:
\begin{itemize}
\item If $\Pi$ is in rows 1--8 of Table~\ref{tab:sixteen}, then there is a drawing $D$ of $K_4$ whose rotation system is $\Pi$, that is, $\Pi$ is realizable. In this case every such drawing $D$ satisfies the following.
\begin{itemize}
\item If $\Pi$ is in rows 1 or 2, then $D$ is {\em plane}, that is, it has no crossings.
\item If $\Pi$ is in rows 3--8, then $D$ has exactly one crossing.
\end{itemize}
\item If $\Pi$ is in rows 9--16, then there is no drawing of $K_4$ whose rotation system is $\Pi$. That is, $\Pi$ is non-realizable.
\end{itemize}
\end{observation}

\begin{proof}[Proof of Observation~\ref{obs:k4rows}]
This is well known, and easy to verify by hand. See for instance~\cite[Observation~14, Table~1 and Figure~4]{kyncl13}, and also~\cite[Observation~2.7]{Sullivan}.
\end{proof}

In view of Observation~\ref{obs:k4rows}, the $16$ rotation systems of size $4$ are naturally classified as follows.

\begin{defbox}{Plane, crossing, and non-realizable rotation systems of size $4$}
Let $Q=\{i,j,k,\ell\}$ be a quad, with $i<j<k<\ell$, so that the rows of Table~\ref{tab:sixteen} are the sixteen rotation systems on $Q$. We say that a rotation system on $Q$ is
\begin{itemize}\setlength{\topsep}{2pt}\setlength{\itemsep}{0pt}
\item {\em plane} if it is in rows 1 or 2 of Table~\ref{tab:sixteen}, and
\item {\em crossing} if it is in rows 3--8.
\end{itemize}
By Observation~\ref{obs:k4rows}, the rotation systems in rows 9--16 are exactly the non-realizable rotation systems on $Q$.
\end{defbox}

Since a quad $Q$ is an ordinary set, in principle it makes no sense to classify $Q$ as plane, crossing, or non-realizable. On the other hand, if we have a rotation system $\Pi$ then we may classify $Q$ according to the type of the rotation subsystem $\Pi|_{Q}$ that $\Pi$ induces on $Q$. This slight abuse of terminology will be very convenient for the rest of this paper.

\begin{defbox}{Plane, crossing, and non-realizable quads (a convenient abuse of terminology)}
Let $\Pi$ be a rotation system of $K_n$. We say that a quad $Q$ is {\em plane}, {\em crossing} or {\em non-realizable in $\Pi$} if the rotation subsystem $\Pi|_Q$ is plane, crossing or non-realizable, respectively.

\medskip
We note that every quad is either plane, crossing or non-realizable in $\Pi$.
\end{defbox}

\subsection{The main result}\label{sub:mainresult}
The main result in this paper involves the following two quantities.

\begin{defbox}{$\crq{\Pi}$ and $\nr(\Pi)$}
Let $\Pi$ be a rotation system of $K_n$. We let $\crq{\Pi}$ be the number of crossing quads in $\Pi$, and $\nr(\Pi)$ the number of non-realizable quads in $\Pi$.
\end{defbox}

 By Observation~\ref{obs:k4rows}, these are the two numbers in the statement of Theorem~\ref{thm:odd} in Section~\ref{sub:ourmain}, which we reproduce for convenience.

\theoremstyle{plain}
\newtheorem*{restatedodd}{Theorem~\ref{thm:odd}}
\begin{restatedodd}[Implies Conjecture~\ref{con:hill}]
Let $\Pi$ be any rotation system of $K_n$, where $n\ge 3$ is an integer. Then
$$
\crq{\Pi}+3\,\nr(\Pi)\ge Z(n).
$$
\end{restatedodd}

We show below that Theorem~\ref{thm:odd} easily implies the Harary--Hill conjecture. The link between drawings and rotation systems is the following remark.

\begin{remark}\label{rem:crD}
Let $D$ be a drawing of $K_n$, and let $\Pi$ be its rotation system. Then $\nr(\Pi)=0$, and $\ucr{D}=\crq{\Pi}$.
\end{remark}

\begin{proof}[Proof of Remark~\ref{rem:crD}]
Let $Q$ be a quad. The drawing of $K_4$ induced by $D$ on $Q$ has rotation system $\Pi|_Q$. By Observation~\ref{obs:k4rows} it follows that $Q$ is not non-realizable in $\Pi$, and $Q$ is crossing in $\Pi$ exactly when this drawing of $K_4$ has a crossing, in which case it has exactly one crossing. Since $D$ is simple, each crossing of $D$ involves two disjoint edges, and so it lies in exactly one of these drawings of $K_4$, namely the one induced on the four endvertices of the two edges. Therefore $\nr(\Pi)=0$, and the crossings of $D$ correspond bijectively to the crossing quads of $\Pi$.
\end{proof}

\begin{corollary}[The Harary--Hill conjecture]\label{cor:hh}
For every integer $n\ge 3$, $\ucr{K_n}=Z(n)$.
\end{corollary}

\begin{proof}[Proof of Corollary~\ref{cor:hh}]
The existence of drawings of $K_n$ with exactly $Z(n)$ crossings implies that $\ucr{K_n}\le Z(n)$, and so it suffices to show that $\ucr{D}\ge Z(n)$ for every drawing $D$ of $K_n$.

Let $D$ be a drawing of $K_n$. By Remark~\ref{rem:crD} the rotation system $\Pi$ of $D$ satisfies $\nr(\Pi)=0$ and $\ucr{D}=\crq{\Pi}$. Since $\nr(\Pi)=0$, Theorem~\ref{thm:odd} implies that $\crq{\Pi}\ge Z(n)$, and so $\ucr{D}\ge Z(n)$.
\end{proof}


\subsection{Settings for the rest of this paper}\label{sub:goalsetting}

The rest of this paper is devoted to the proof of Theorem~\ref{thm:odd}. 

\vglue 0.3 cm
\noindent{\bf Remark. }{\em It suffices to prove Theorem~\ref{thm:odd} for odd $n$. Indeed, the case of even $n$ follows from the case of odd $n$ by the same counting argument that proves Observation~\ref{obs:oddn} (see for instance~\cite{szekely}).}
\vglue 0.3 cm

In view of this remark, for the rest of the paper we assume that $n$ is odd. For $n=3$ there is nothing to prove, since $Z(3)=0$. For $n\ge 5$ we work under the following setting.

\begin{settingsbox}{Settings for the rest of this paper}
Throughout this paper:
\par\vspace{8pt}
\begin{itemize}\setlength{\itemsep}{6pt}
\item $\Pi=\bigl(\pi(1),\ldots,\pi(n)\bigr)$ is a fixed rotation system of $K_n$, where $n\ge 5$ is odd.
\item $r:=(n-1)/2$.
\end{itemize}
\end{settingsbox}


\section{Proof of Theorem~\ref{thm:odd} from 35,000 feet (Optional)}\label{sub:sketch}

Our aim in this section is to discuss the strategy, main steps, and ingredients in the proof of Theorem~\ref{thm:odd}. We emphasize that the rest of the paper is independent of this section, and so the reader who wants to move on right away to the proof may safely skip this section.

As we are about to see, much of the proof of Theorem~\ref{thm:odd} revolves around facts (some well-known, some established in~\cite{openai}) about spaces of polynomials. Now a combinatorialist interested in drawings and crossing numbers will likely have the same reaction we had when we started reading~\cite{openai}: why are polynomials relevant to crossing numbers, and in particular to the Harary--Hill conjecture? We address this question in Section~\ref{sub:whypoly}.


\subsection{The strategy: a linear map and a dimension count}\label{sub:strategy}

At a very high level the strategy of the proof is as follows. 

The proof involves a set $\xx$ of ordered pairs $(ij,k\ell)$ of disjoint edges, determined by $\Pi$, of size $|\xx|=2\crq{\Pi}+6\,\nr(\Pi)$. Its precise definition is given in Section~\ref{sec:keyobj}, and is not needed at all for this discussion. 

We define the {\em evaluation map} 
$$E:\ss\to \real^{\xx},$$ 
a linear map whose domain is a space $\ss$ of polynomials, and that sends a polynomial to a vector in $\real^{\xx}$.

More precisely, this map sends a polynomial $P=P(x_1,x_2,y_1,y_2)\in\ss$ to the list of its values $P(i,j,k,\ell)$ at the ordered pairs $(ij,k\ell)\in\xx$, and ignores all the other values of $P$.

The map $E$ is how the specific rotation system $\Pi$ under consideration enters the picture. As we shall see, the space $\ss$ of polynomials depends only on $n$, but $E$ involves the crossing quads and the non-realizable quads of $\Pi$, and so it depends on $\Pi$.

We remark that since $E$ depends on $\Pi$, in principle it would be appropriate to emphasize this by writing $E_\Pi$. But this would put an unnecessary load on the notation, as we will always be working with an arbitrary but fixed rotation system $\Pi$.

It is easily proved that the dimension of $\ss$ is $Z(n)$, that is, (i) $\dim \ss=Z(n)$. Theorem~\ref{thm:odd} follows from this and two other facts, namely that (ii) $E$ is injective (that is, a polynomial in $\ss$ that vanishes at every element of $\xx$ must be the zero polynomial) and that (iii) the image of $\ss$ under $E$ has dimension at most $\frac{1}{2}|\xx|$. 

To see that indeed Theorem~\ref{thm:odd} follows from these three facts, first note that (ii) implies that $\dim E(\ss)=\dim\ss$ and (iii) implies that $\frac{1}{2}|\xx|\ge \dim E(\ss)$, that is, $\crq{\Pi}{}+3\,\nr(\Pi) \ge \dim E(\ss)$. Putting together (i), (ii), and (iii) it follows that $\crq{\Pi}{}+3\,\nr(\Pi) \ge \dim E(\ss)=\dim\ss = Z(n)$, as claimed in Theorem~\ref{thm:odd}. 

Facts (i), (ii), and (iii) are Lemma~\ref{lem:dimS} and Propositions~\ref{pro:inj} and~\ref{pro:iso}, respectively. These are formally stated in Section~\ref{sec:reduction}, where for completeness we show that Theorem~\ref{thm:odd} follows from them.


\subsection{The workhorse: the Main Identity}\label{sub:workhorse}

The workhorse behind the proofs of both propositions is a single identity, which we call the {\em Main Identity}. This identity, introduced as Theorem~\ref{thm:master} in Section~\ref{sec:masterid}, relates the crossing and non-realizable quads of $\Pi$ to the cyclic permutations $\pi(1),\ldots,\pi(n)$ that constitute $\Pi$.

The Main Identity involves a polynomial $A$ in four variables, and it holds for every polynomial $A$ of a certain kind. So it is not a single identity but a whole family of them, one for each $A$. As we will see in Sections~\ref{app:inj} and~\ref{sec:iso}, each of the two propositions is proved by choosing $A$ cleverly, so that one side of the identity vanishes and the identity becomes a statement about the other side.


\subsection{Why polynomials?}\label{sub:whypoly}

This (at first sight inscrutable) Main Identity is an identity written in terms of polynomials. Why polynomials at all? We devote this subsection to a mostly informal discussion around this natural question.

The simplest way to count the crossing quads of $\Pi$ is to let each of them contribute $1$. In the Main Identity, a crossing quad does not simply contribute $1$. Instead, it contributes a number that depends on its four vertices, and that may be different for different crossing quads.

A convenient way to choose such numbers is to take a polynomial $A=A(x_1,x_2,y_1,y_2)$ in four variables. Each crossing quad has one distinguished pair of disjoint edges $ij$ and $k\ell$ (in the language of Section~\ref{sub:quads}, the only one with nonzero sign). Since the vertices are the numbers $1,\ldots,n$, we can let the quad contribute $A(i,j,k,\ell)$, multiplied by a sign and a fixed weight. The non-realizable quads contribute in a similar way.

The rotations enter in the same way. Each pair of edges $ij$ and $ik$ leaving the same vertex $i$ also contributes a number, namely $A(i,j,i,k)$, again multiplied by a sign and a fixed weight. In a nutshell, the Main Identity says that the total contribution of the quads and the total contribution of the pairs of edges at the vertices add up to zero.

It seems worth pointing out a similarity with the discharging method, a powerful technique frequently used in topological graph theory. In a discharging argument every vertex and face receives a number (namely its {\em charge}), and the charges are then moved around by rules that at first sight look arbitrary. What makes the argument work is a fixed identity, Euler's formula, which says what the total charge must be, however the charges are redistributed.

The Main Identity plays a similar role here. The numbers contributed by the quads and by the pairs of edges at a vertex are the charges, the polynomial $A$ is the rule that assigns them, and the identity says that the two totals always add up to zero. As in a discharging argument, the choice of the rule, that is, of $A$, looks arbitrary until one sees the powerful consequences of this choice.


\subsection{The proof of the Main Identity}\label{sub:topology}

After Section~\ref{sec:masterid} the rest of the paper is devoted to the proofs of the Main Identity (Theorem~\ref{thm:master}) and of Propositions~\ref{pro:inj} and~\ref{pro:iso}.

The Main Identity is proved in Sections~\ref{sec:proofmainid} and~\ref{sec:polysec}. Most of the work is done in Section~\ref{sec:proofmainid}, where we establish a related identity (the {\em Cycle Identity}, stated as Theorem~\ref{thm:cancel}) that relates the signs of the pairs of disjoint edges formed by two closed walks to the rotations at their common vertices. As we show in Section~\ref{sec:polysec}, the Main Identity is a rather straightforward consequence of the Cycle Identity. One could actually argue that the Main Identity is ``the polynomial version'' of the Cycle Identity.


\subsection{The proofs of the two propositions}\label{sub:twoprops}

The proofs of Propositions~\ref{pro:inj} and~\ref{pro:iso} are in Sections~\ref{app:inj} and~\ref{sec:iso}. The core of each of these proofs is an application of the Main Identity to a suitably chosen polynomial $A$. 

The proof of Proposition~\ref{pro:inj} has an additional major ingredient, namely Lemma~\ref{lem:detector}. This lemma is a purely algebraic statement about polynomials in two variables. The proof of the lemma is the content of Section~\ref{sec:interp}.


\section{Reducing the proof of Theorem~\ref{thm:odd} to two propositions}\label{sec:reduction}

The proof of Theorem~\ref{thm:odd} revolves around a linear map, the {\em evaluation map} $E$. In this section we define this map, state one lemma and two propositions around it, and show that together they easily imply Theorem~\ref{thm:odd}. We prove the lemma at the end of the section. The rest of the paper is devoted to the proofs of the two propositions.


\subsection{Directed edges, orientations, and the signs of a quad}\label{sub:quads}

The proof uses a lot more information about a quad than its type. For each of the three pairs of disjoint edges of the quad, it uses a number that we call the {\em sign} of the pair. 

In order to define this number we first need to lay out the terminology we use for directed edges, as well as the notion of the orientation of a vertex in a quad.

\begin{defbox}{Directed edges}
We write $ij$ for the edge joining $i$ and $j$, directed from $i$ to $j$:\quad
\raisebox{-0.55cm}{\includegraphics{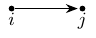}}\\
Thus $ij$ and $ji$ are the same edge, with opposite directions.
\end{defbox}

In the next definition we use cyclic orders. If $Q$ is a quad and $i$ is a vertex of $Q$, then the other three vertices of $Q$ have two cyclic orders. In one of these cyclic orders, if we start at the smallest element the three elements appear in increasing order. We call this cyclic order {\em increasing}, and we call the other cyclic order {\em decreasing}.

\begin{defbox}{The orientation of a vertex in a quad}
Let $Q$ be a quad and let $i$ be a vertex of $Q$. The {\em orientation of $i$ in $Q$} is $\orient{i}(Q)=+1$ if the rotation at $i$ in $\Pi|_Q$ is the increasing cyclic order, and $\orient{i}(Q)=-1$ if it is the decreasing one. When $Q$ is clear from the context we simply write $\orient{i}$.
\end{defbox}

For instance, in row 1 of Table~\ref{tab:sixteen} the rotations at $i$ and $k$ are increasing and those at $j$ and $\ell$ are decreasing, so $(\orient{i},\orient{j},\orient{k},\orient{\ell})=(+1,-1,+1,-1)$. In row 9 the rotations at $i$, $j$ and $k$ are increasing and the one at $\ell$ is decreasing, so $(\orient{i},\orient{j},\orient{k},\orient{\ell})=(+1,+1,+1,-1)$.

We emphasize that the orientation of a vertex depends on the quad, and not only on the vertex. A vertex lies in many quads, and it may have orientation $+1$ in some of them and $-1$ in others.

For instance, suppose that $\pi(1)=(2,5,4,7,6,3)$ as in the illustration in Section~\ref{sub:rotsys}. In the quad $\{1,2,3,4\}$ the rotation at $1$ is $(2,4,3)$, which is decreasing, and so the orientation of $1$ in this quad is $-1$. In the quad $\{1,2,4,7\}$ the rotation at $1$ is $(2,4,7)$, which is increasing, and so the orientation of $1$ in this quad is $+1$.

\begin{defbox}{The three signs of a quad}
Let $Q=\{i,j,k,\ell\}$ be a quad, with $i<j<k<\ell$, and let $\orient{i},\orient{j},\orient{k},\orient{\ell}$ be the orientations of its four vertices. The {\em signs} of the three pairs of disjoint edges of $Q$ are
\begin{equation}\label{eq:signs}
\begin{aligned}
\sign(ij,k\ell)&{}=\tfrac14\big(-\orient{i}+\orient{j}+\orient{k}-\orient{\ell}\big),\\
\sign(ik,j\ell)&{}=\tfrac14\big(\orient{i}+\orient{j}+\orient{k}+\orient{\ell}\big),\\
\sign(i\ell,jk)&{}=\tfrac14\big(-\orient{i}-\orient{j}+\orient{k}+\orient{\ell}\big).
\end{aligned}
\end{equation}
\end{defbox}

\begin{figure}[h!]
\centering
\includegraphics{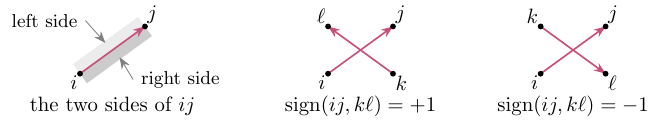}
\caption{The sign of a crossing in a drawing. If the directed edges $ij$ and $k\ell$ cross, the sign of the ordered pair $(ij,k\ell)$ is $+1$ if $k\ell$ crosses $ij$ from its right side, and $-1$ if $k\ell$ crosses $ij$ from its left side.}
\label{fig:sign}
\end{figure}

It seems worth explaining where this definition comes from. Suppose that $\Pi|_Q$ is realizable, that is, it is the rotation system of a drawing of $K_4$. In such a drawing either no two edges cross, or exactly two disjoint edges cross, and in this case the crossing has a sign that tells from which side one edge crosses the other edge (see Figure~\ref{fig:sign}). The key observation is that whether or not two edges cross, and the sign of their crossing if it exists, are determined by the rotations at the four vertices (see~\cite[Lemma~10]{kyncl13}, and also~\cite[Observation~2.7]{Sullivan}), and so they are determined by the four orientations $\orient{i}$, $\orient{j}$, $\orient{k}$, and $\orient{\ell}$.

The formulas above are the {\em only} linear combinations of $\orient{i}$, $\orient{j}$, $\orient{k}$, and $\orient{\ell}$ that give the right signs on realizable quads (that is, they give zero if a pair of edges does not cross each other, and the sign of the crossing if the edges do cross). On a non-realizable quad there is no associated drawing, and so there is no crossing whose sign we could record. We use the same formulas anyway, and this is not just a matter of consistency. As it happens, the rotations at the four vertices of a quad do not allow any other choice, since the identity behind the whole proof (namely the Cycle Identity, Theorem~\ref{thm:cancel}) holds with these numbers as ``signs'' and fails with any others.

\vglue 0.3 cm
\noindent{\bf Remark. }{\em The previous discussion, including Figure~\ref{fig:sign}, only motivates the definition of the signs. Neither the definition nor the proof of Theorem~\ref{thm:odd} depends on drawings.}
\vglue 0.3 cm

\begin{defbox}{Convention: signs of ordered pairs of directed edges}
 The previous box gives the sign of each pair $\{ij,k\ell\}$ of disjoint edges of $Q$ with $i<j$, $k<\ell$ and $i<k$, as the sign $\sign(ij,k\ell)$ of the ordered pair of directed edges $(ij,k\ell)$. For any other order of the two edges, or any other direction of either of them, the sign is given by the rules
$$
\sign(k\ell,ij)=\sign(ji,k\ell)=\sign(ij,\ell k)=-\sign(ij,k\ell).
$$
That is, exchanging the two edges, or reversing either of them, changes the sign.
\end{defbox}

We make essential use of the fact that the type of a quad can be recognized from its three signs. We formalize this in the next statement. 

\begin{remark}\label{rem:signsbytype}
A plane quad has all three signs equal to $0$. A crossing quad has exactly one nonzero sign, which is either $+1$ or $-1$. Finally, a non-realizable quad has three signs, each of which is either $+\tfrac12$ or $-\tfrac12$.
\end{remark}

\begin{proof}[Proof of Remark~\ref{rem:signsbytype}]
We use the following simple fact. Consider any four numbers, each of which is either $+1$ or $-1$. If an odd number of these numbers are $-1$, then their sum is $\pm 2$, and if an even number of them are $-1$, then their sum is $0$ or $\pm 4$.

Each of the three sums in~\eqref{eq:signs} adds up the four orientations, with either none or two of them multiplied by $-1$. Multiplying two of the terms by $-1$ does not change whether an even or an odd number of terms are $-1$. A glance at Table~\ref{tab:sixteen} shows that an even number of the four orientations are $-1$ in rows 1--8, and an odd number in rows 9--16.

In rows 1--8 each of the three sums is therefore $0$ or $\pm 4$, and so each sign is $0$ or $\pm 1$. Moreover, it is easily verified that all three signs are $0$ in rows 1 and 2, and that exactly one of them is nonzero in rows 3--8. This proves the remark for the plane and the crossing quads.

Finally, in rows 9--16 each sum is $\pm 2$, and so each sign is $\pm\tfrac12$. Thus the remark also holds for non-realizable quads.
\end{proof}


\subsection{The key object: the linear map $E$}\label{sec:keyobj}

We start by introducing the two ingredients of $E$, the set $\xx$ and the space $\ss$.

We recall from Remark~\ref{rem:signsbytype} that each of the $\crq{\Pi}$ crossing quads has exactly one pair of disjoint edges with nonzero sign, and this sign is $\pm 1$. Each of the $\nr(\Pi)$ non-realizable quads has three pairs of disjoint edges, with signs $\pm\tfrac12$. All other pairs of disjoint edges have sign $0$.

\begin{defbox}{The set $\xx$ of pairs with nonzero sign, each pair appearing twice}
We let $\xx$ be the set of ordered pairs $(ij,k\ell)$ of disjoint edges such that $\sign(ij,k\ell)\ne 0$, with the convention that each edge is written with its smaller endpoint first, that is, $i<j$ and $k<\ell$. 
\end{defbox}

Since each pair with nonzero sign involves two ordered pairs, namely $(ij,k\ell)$ and $(k\ell,ij)$, the set $\xx$ is twice as large as the number of such pairs, which is $\crq{\Pi}+3\,\nr(\Pi)$. We record this remark for later use.

\begin{remark}\label{rem:twice}
We have $|\xx|=2\,\crq{\Pi}{}+6\,\nr(\Pi)$.
\end{remark}

\label{sec:pairs}
Polynomials play a central role in the proof of Theorem~\ref{thm:odd}. The key object is defined on a space of polynomials, and the workhorse behind the proof is an identity that holds for every polynomial of a certain kind. We now give the relevant definitions.

\begin{defbox}{Pair polynomials}
A {\em pair polynomial} is a polynomial $P(x_1,x_2,y_1,y_2)$ in four variables, of degree at most $n-2$ in each variable.

\medskip
Since the vertices are the numbers $1,\ldots,n$, we can evaluate a pair polynomial at any ordered pair of edges $(ij,k\ell)$, and say that the value of $P$ at the pair $(ij,k\ell)$ is simply $P(i,j,k,\ell)$. This explains why we call them ``pair'' polynomials.
\end{defbox}

 For example, let $n\ge 5$ and $P=x_1^2+x_2y_1-3y_2$. Thus $P$ is a pair polynomial, as the degree in each variable is not greater than $5-2=3$. The value of $P$ at the ordered pair $(15,24)$ of edges is $P(1,5,2,4)=1+10-12=-1$, and at the ordered pair $(24,15)$ we have $P(2,4,1,5)=4+4-15=-7$.

\begin{defbox}{Symmetric pair polynomials}
\setlength{\abovedisplayskip}{12pt}\setlength{\belowdisplayskip}{12pt}
A pair polynomial $P$ is {\em symmetric} if its value is not affected by the exchanges $x_1 \leftrightarrow x_2$ and/or $y_1\leftrightarrow y_2$. That is, $P$ is symmetric if
$$
P(x_2,x_1,y_1,y_2)=P(x_1,x_2,y_1,y_2)\quad\text{and}\quad P(x_1,x_2,y_2,y_1)=P(x_1,x_2,y_1,y_2).
$$
\end{defbox}

 For example, $x_1+x_2$ and $x_1x_2(y_1+y_2)^2$ are symmetric, and $x_1^2+x_2y_1-3y_2$ is not.

\begin{defbox}{The space $\ss$}
We let $\ss$ be the space of symmetric pair polynomials of degree at most $r-2$ in each variable.
\end{defbox}

 For example, for $n=5$ we have $r=2$. Thus in this case the degree in each variable is at most $0$, and so the space $\ss$ consists of the constant polynomials and $\dim\ss=1$. For $n=7$ we have $r=3$, and so the degree in each variable is at most $1$. The polynomials of degree at most $1$ in $x_1$ and in $x_2$ that are symmetric in $x_1,x_2$ are the combinations of
$$
1,\qquad x_1+x_2,\qquad x_1x_2,
$$
and the same three expressions in $y_1,y_2$ for the second pair. Thus in this case $\ss$ has as a basis the $9$ products of one expression from each list, such as $(x_1+x_2)\,y_1y_2$, and $\dim\ss=9$.

Having introduced $\ss$ and $\xx$, we can now define the key object in the proof of Theorem~\ref{thm:odd}.

\begin{defbox}{The key object: the linear map $E$}
The {\em evaluation map} is the linear map
$$
\scalebox{1.2}{$E:\ss\longrightarrow\real^{\xx}$}
$$
that sends each $P\in\ss$ to the vector of its values at the elements of $\xx$. That is,
$$
E(P)=\big(P(i,j,k,\ell)\big)_{(ij,k\ell)\in\xx}\quad.
$$
\end{defbox}
\medskip


\subsection{Reducing Theorem~\ref{thm:odd} to two propositions}\label{sec:workhorses}

We now state the lemma and the two propositions mentioned above.

We start by stating the lemma, which unveils the connection between $\ss$ and $Z(n)$. As usual, we write $\dim V$ for the dimension of a vector space $V$.

\begin{lemma}\label{lem:dimS}
$\dim\ss=Z(n)$.
\end{lemma}

Lemma~\ref{lem:dimS} is the first and last place in which the number $Z(n)$ makes an appearance in the proof of Theorem~\ref{thm:odd}. The long and difficult parts of the argument are the proofs of the next two propositions, and $Z(n)$ never shows up in any of them.

The next two propositions are the heart of the proof of Theorem~\ref{thm:odd}.

\begin{proposition}\label{pro:inj}
The evaluation map $E$ is injective. That is, if $P\in\ss$ satisfies that $P(i,j,k,\ell)=0$ for every pair $(ij,k\ell)\in\xx$, then $P$ is the zero polynomial.
\end{proposition}

\begin{proposition}\label{pro:iso}
$\displaystyle\dim E(\ss)\le\tfrac{1}{2}|\xx|.$ 
\end{proposition}

For completeness let us formally show that together with Lemma~\ref{lem:dimS} these propositions imply Theorem~\ref{thm:odd}.

\begin{proof}[Proof of Theorem~\ref{thm:odd}, assuming Propositions~\ref{pro:inj} and~\ref{pro:iso}]
We recall from Section~\ref{sub:goalsetting} that the case $n=3$ is trivial, and that the even $n$ case follows from the odd $n$ case. Thus we may assume that $n\ge 5$ is odd, as in the settings we have fixed since Section~\ref{sub:goalsetting}.

The map $E$ is linear, and so its image $E(\ss)$ is a vector space. By Proposition~\ref{pro:inj} the map $E$ is injective, and so $\dim E(\ss)=\dim\ss$. Now by Lemma~\ref{lem:dimS} and Proposition~\ref{pro:iso} we have $\dim\ss=Z(n)$ and  $\frac{1}{2} |\xx|\ge \dim E(\ss)$, respectively. Putting these three facts together and using Remark~\ref{rem:twice} we obtain
$$
\crq{\Pi}+3\,\nr(\Pi)=\tfrac{1}{2}|\xx|\ge\dim E(\ss)=\dim\ss=Z(n).\qedhere
$$
\end{proof}

\begin{proof}[Proof of Lemma~\ref{lem:dimS}]
Let $P=P(x_1,x_2,y_1,y_2)\in\ss$, and let $c$ be the coefficient in $P$ of a monomial $x_1^{s_1}x_2^{s_2}y_1^{t_1}y_2^{t_2}$, where all four exponents are between $0$ and $r-2$. Since $P$ is symmetric, the monomials obtained from $x_1^{s_1}x_2^{s_2}y_1^{t_1}y_2^{t_2}$ by exchanging $x_1$ and $x_2$, or $y_1$ and $y_2$, or both, also have coefficient $c$ in $P$. Thus $P$ is determined by its coefficients of the monomials with $s_1\le s_2$ and $t_1\le t_2$. Conversely, every choice of these coefficients gives a polynomial in $\ss$. 

From this it follows that $\dim\ss$ is the number of monomials with $0\le s_1\le s_2\le r-2$ and $0\le t_1\le t_2\le r-2$. An elementary counting argument shows that there are $\binom{r}{2}$ pairs $(s_1,s_2)$ such that $0\le s_1\le s_2\le r-2$, and similarly $\binom{r}{2}$ pairs $(t_1,t_2)$ such that $0\le t_1\le t_2\le r-2$. Therefore
$$
\dim\ss=\binom{r}{2}^2=\frac{r^2(r-1)^2}{4}=Z(2r+1).\qedhere
$$
\end{proof}


\section{The Main Identity}\label{sec:masterid}

In this section we state the Main Identity, the workhorse behind the proofs of Propositions~\ref{pro:inj} and~\ref{pro:iso}, and sketch how we use the identity in the proofs.

The Main Identity involves a few more basic concepts, namely the rotation signs and the weights of the vertices. Let us start by introducing these elementary notions.


\subsection{Rotations, signs and weights}\label{sec:conv}

We recall that rotations are written as linear permutations starting at their smallest element (Section~\ref{sub:rotsys}). We also recall that the sign $\sign(ij,k\ell)$ of two disjoint edges $ij$ and $k\ell$ is defined in Section~\ref{sub:quads}.

For convenience we start by extending the sign to pairs of edges that share an endvertex.

\begin{defbox}{The sign of two edges that share an endvertex}
If $ij$ and $k\ell$ share an endvertex, that is, if $\{i,j\}\cap \{k,\ell\}\neq\emptyset$, then we let $\sign(ij,k\ell)=0$. Thus for all $i<j$ and $k<\ell$, we have $\sign(ij,k\ell)\ne 0$ exactly when $(ij,k\ell)\in\xx$.
\end{defbox}

We now define a crucial notion, the rotation sign. It records whether one vertex comes before or after another in the rotation at a third vertex.

\begin{defbox}{The rotation sign $\tau_i(j,k)$}
Let $i$ be a vertex, and let $j\neq k$ be vertices different from $i$. The {\em rotation sign} $\tau_i(j,k)$ is $+1$ if $j$ occurs before $k$ in $\pi(i)$, and $-1$ if $j$ occurs after $k$. For convenience we also let $\tau_i(j,j)=0$ for every vertex $j\ne i$.
\end{defbox}

\begin{defbox}{The weight $\omega_i$ of vertex $i$}
\setlength{\abovedisplayskip}{14pt}\setlength{\belowdisplayskip}{14pt}
\setlength{\abovedisplayshortskip}{14pt}\setlength{\belowdisplayshortskip}{14pt}
The weight $\omega_i$ of vertex $i$ is
$$
\omega_i=\frac{1}{\prod_{j\in [n]\setminus\{i\}}(i-j)}.
$$
\end{defbox}

 For example, for $n=5$ the weights are $\omega_1=\tfrac{1}{24}$, $\omega_2=-\tfrac{1}{6}$, $\omega_3=\tfrac{1}{4}$, $\omega_4=-\tfrac{1}{6}$ and $\omega_5=\tfrac{1}{24}$.


\subsection{The workhorse: the Main Identity}\label{sec:mainid}
We need one last notion in order to state the Main Identity. The identity holds for pair polynomials of a particular kind, namely those that change sign when either pair of variables is exchanged. Formally we have the following definition.

\begin{defbox}{Antisymmetric pair polynomials}
\setlength{\abovedisplayskip}{12pt}\setlength{\belowdisplayskip}{12pt}
A pair polynomial $P=P(x_1,x_2,y_1,y_2)$ is {\em antisymmetric} if it changes sign under each of the exchanges $x_1 \leftrightarrow x_2$ and $y_1\leftrightarrow y_2$. That is, $P$ is antisymmetric if
$$
P(x_2,x_1,y_1,y_2)=-P(x_1,x_2,y_1,y_2)\quad\text{and}\quad P(x_1,x_2,y_2,y_1)=-P(x_1,x_2,y_1,y_2).
$$
\end{defbox}

 For example, the polynomial $(x_1-x_2)(y_1-y_2)(x_1+x_2)$ is antisymmetric.

\begin{theorem}[Main Identity]\label{thm:master}
For every antisymmetric pair polynomial $A$,
\begin{equation}\label{eq:master}
\vcenter{\hbox{\includegraphics{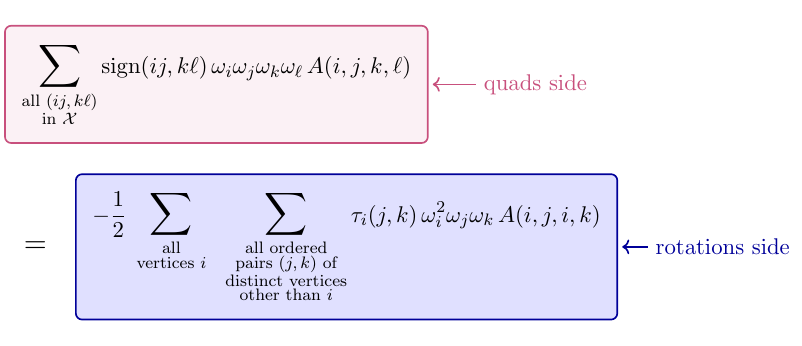}}}
\end{equation}
\end{theorem}
\medskip

The proof of Theorem~\ref{thm:master} takes roughly half the rest of this paper. In Section~\ref{sec:proofmainid} we state and prove a related identity (the {\em Cycle Identity}, see Theorem~\ref{thm:cancel}), and in Section~\ref{sec:polysec} we show that the Main Identity is a rather straightforward consequence of the Cycle Identity. (For the reader who wants to compare, the Main Identity and the Cycle Identity are versions of Lemmas~4.1 and~2.1 in~\cite{openai}, respectively.)


\section[Proof of the Main Identity I: a simpler version without polynomials]{Proof of the Main Identity I: a simpler version\\without polynomials}\label{sec:proofmainid}

The Main Identity of Section~\ref{sec:mainid} is the polynomial version of a simpler identity about two closed walks in $K_n$, which involves no polynomials at all. We call this the {\em Cycle Identity} (as opposed to the longer but plausible ``Closed Walks Identity'') for reasons that we explain in Remark~\ref{rem:cycles}. In this section we state and prove this identity, and in the next section we use it to prove the Main Identity.

We recall that since the end of Section~\ref{sub:goalsetting} we work with an arbitrary but fixed rotation system $\Pi$ of $K_n$. The two ingredients of the Cycle Identity are the signs of Section~\ref{sub:quads} and the rotation signs $\tau_i$ of Section~\ref{sec:conv}.  

\subsection{Encoding a closed walk into a matrix: balanced matrices}\label{sec:cycid}

The Cycle Identity is a theorem (namely Theorem~\ref{thm:cancel}) about closed walks in $K_n$. Our first step towards stating this theorem is to describe a closed walk in a way that is suitable for computations.

We proceed as follows. For every ordered pair $(i,j)$ of distinct vertices, we record how many times the walk goes from $i$ to $j$ minus how many times it goes from $j$ to $i$. This gives a table of numbers $C(i,j)$, that is, an $n\times n$ matrix. Figure~\ref{fig:walkmatrices} shows two examples for $n=5$.

\begin{figure}[h!]
\centering
$
\begin{array}{c|ccccc}
 & 1 & 2 & 3 & 4 & 5\\ \hline
1 & 0 & 1 & 0 & -1 & 0\\
2 & -1 & 0 & 0 & 1 & 0\\
3 & 0 & 0 & 0 & 0 & 0\\
4 & 1 & -1 & 0 & 0 & 0\\
5 & 0 & 0 & 0 & 0 & 0
\end{array}
\qquad\qquad
\begin{array}{c|ccccc}
 & 1 & 2 & 3 & 4 & 5\\ \hline
1 & 0 & 1 & 0 & -1 & 0\\
2 & -1 & 0 & 1 & 0 & 0\\
3 & 0 & -1 & 0 & 1 & 0\\
4 & 1 & 0 & -1 & 0 & 0\\
5 & 0 & 0 & 0 & 0 & 0
\end{array}
$
\caption{The matrices of two closed walks in $K_5$. On the left-hand side we have the matrix of the closed walk (a triangle) $1\to 2\to 4\to 1$. On the right-hand side we have the matrix of the closed walk $1\to 2\to 3\to 1\to 3\to 4\to 1$. This walk goes once from $3$ to $1$ and once from $1$ to $3$, and so the entries at $(1,3)$ and $(3,1)$ are $0$.}
\label{fig:walkmatrices}
\end{figure}

Two properties of these matrices are immediately clear. First, $C(j,i)=-C(i,j)$. Second, the entries of every row add up to zero, since a closed walk leaves each vertex as many times as it arrives at the vertex. Working with matrices has the advantage that a walk that passes through a vertex several times does not need any special treatment. An analogous remark applies for two walks that share an edge.

An important additional advantage of working with matrices is that they allow arbitrary real entries, which are needed in Section~\ref{sec:polysec}. There the Cycle Identity is applied to matrices that do not come from any closed walk.

The following notion captures the properties shared by all matrices obtained from closed walks in the way we have just described.

\begin{defbox}{Balanced matrices}
Let $C$ be a real $n\times n$ matrix whose rows and columns are indexed by the vertices. We say that $C$ is {\em balanced} if it is antisymmetric (that is, $C(j,i)=-C(i,j)$ for all $i,j\in[n]$) and the sum of the entries of every row is zero, that is,
\begin{equation}\label{eq:balance}
\textstyle\sum_j C(i,j)=0\qquad\text{for every }i \in [n].
\end{equation}
It seems natural to regard $C(i,j)$ as an amount flowing along the edge from $i$ to $j$. Balancedness then means that at every vertex the amount of incoming flow is the same as the amount of outgoing flow, just as a closed walk arrives at each vertex as many times as it leaves it.
\end{defbox}

\begin{remark}[Balanced matrices and cycles in algebraic topology]\label{rem:cycles}
In the language of algebraic topology (see for instance~\cite{giblin}) a balanced matrix corresponds to a {\em $1$-cycle} of the graph, that is, a formal combination of directed edges with zero boundary. Having zero boundary is the balance condition~\eqref{eq:balance} at every vertex. This is where the allusion to cycles in ``Cycle Identity'' comes from. In short, closed walks are cycles (in algebraic topology), and cycles are balanced matrices.
\end{remark}

 The simplest nontrivial closed walks are triangles. As it happens, triangles are also the ``building blocks'' of all closed walks. Every balanced matrix, and in particular every closed walk, can be assembled from triangles. This is why, as we shall see, in order to prove the Cycle Identity it will be enough to focus on triangles.

Matrices of triangles occur so frequently in our discussions that we need a definition for them.

\begin{defbox}{Triangle matrices}
For distinct vertices $i,j,k$, the {\em triangle matrix} $T_{ijk}$ is the $n\times n$ matrix, with rows and columns indexed by the vertices, that has entries $1$ at $(i,j)$, $(j,k)$, $(k,i)$, entries $-1$ at $(j,i)$, $(k,j)$, $(i,k)$, and $0$ elsewhere.
\end{defbox}

It is straightforward to check that $T_{ijk}$ is a balanced matrix, and that it describes the closed walk $i\to j\to k\to i$.

Let us now make precise what we meant when we said that triangles are the ``building blocks'' of all closed walks.

\begin{lemma}\label{lem:span}
Every balanced matrix can be written as a linear combination of the triangle matrices $T_{1jk}$, $2\le j<k\le n$.
\end{lemma}

\begin{proof}[Proof of Lemma~\ref{lem:span}]
Let $C$ be a balanced matrix. In a nutshell, the idea is to subtract from $C$ a suitable linear combination of the matrices $T_{1jk}$ and then show that what is left is the zero matrix.

We show that
$$
C=\sum_{2\le j<k\le n}C(j,k)\,T_{1jk}.
$$
To prove this we let
\begin{equation}\label{eq:Cprime}
C_0:=C-\sum_{2\le j<k\le n}C(j,k)\,T_{1jk},
\end{equation}
and show that $C_0$ is the zero matrix. We note that since $C$ and the triangle matrices are balanced it follows that $C_0$ is also balanced.

We start by showing that ($*$) {\em all the entries of $C_0$ that are not in row $1$ or column $1$ are zero.} To prove this let $j,k$ be any two fixed integers such that $2\le j<k\le n$. Among the triangle matrices in the sum in~\eqref{eq:Cprime} the matrix $T_{1jk}$ is the only one that has a nonzero entry at $(j,k)$, and this entry is $1$. Therefore $C_0(j,k)=C(j,k)-C(j,k)=0$. Since $C_0$ is antisymmetric we also have that $C_0(k,j)=0$, so $C_0(j,k)=0$ for all distinct $j,k \in \{2,\ldots,n\}$. Combining this with the fact that (by antisymmetry) $C_0(j,j)=0$ for all $j\in[n]$ we obtain ($*$).

Now let $j\in[n]$ with $j\ne 1$. In view of ($*$) the only entry in row $j$ of $C_0$ that may be nonzero is $C_0(j,1)$. Since $C_0$ is balanced it follows that the sum of the entries of row $j$ must be zero, and so $C_0(j,1)=0$. Since this holds for every $j\ne 1$ and $C_0(1,1)=0$ (as $C_0$ is antisymmetric), it follows that ($**$) {\em all the entries of column $1$ are zero}.

Using ($**$) and that $C_0$ is antisymmetric we obtain that ($***$) {\em all the entries of row $1$ must also be zero}. Combining ($*$), ($**$) and ($***$) it follows that $C_0$ is the zero matrix, as claimed.
\end{proof}


\subsection{The quad count and the rotation count}\label{sec:cycformal}

 With closed walks written as balanced matrices, the two quantities that appear in the Cycle Identity are  the {\em quad count} $\Cr(C,C')$ and the {\em rotation count} $\Rot(C,C')$, which we now define. 

\begin{defbox}{The quad count and the rotation count}
\setlength{\abovedisplayskip}{14pt}\setlength{\belowdisplayskip}{14pt}
\setlength{\abovedisplayshortskip}{14pt}\setlength{\belowdisplayshortskip}{14pt}
Let $C,C'$ be antisymmetric matrices. We let
\begin{equation}\label{eq:J}
\Cr(C,C')=\sum_{i<j}\,\sum_{k<\ell}C(i,j)\,C'(k,\ell)\,\sign(ij,k\ell)
\end{equation}
and
\begin{equation}\label{eq:Rot}
\Rot(C,C')=\tfrac{1}{2}\sum_i\,\sum_{j,k\ne i}C(i,j)\,C'(i,k)\,\tau_i(j,k).
\end{equation}
 Recall that, for $i<j$ and $k<\ell$, we have $\sign(ij,k\ell)=0$ when $(ij,k\ell)\notin\xx$. Therefore $\Cr$ is really a sum over the elements $(ij,k\ell)$ of the set $\xx$ of Section~\ref{sec:keyobj}, and $\Rot$ is a sum over pairs of edges with a common endvertex.
\end{defbox}

\subsubsection{Intuition behind $\Cr$ and $\Rot$ (optional)}\label{sub:intuition}

Even though the proof of Theorem~\ref{thm:odd} does not use drawings, both~\eqref{eq:J} and~\eqref{eq:Rot} can be traced back to definitions in~\cite{openai} that involve drawings, and it seems worth discussing where they come from. 

Suppose throughout this discussion that $\Pi$ is the rotation system of a drawing $D$ of $K_n$, and that $C$ and $C'$ are the matrices of two closed walks $W$ and $W'$. In the drawing $D$ these walks are (not necessarily simple) closed curves, and it is a well-known consequence of the Jordan curve theorem that their crossings, counted with signs, add up to zero. Formally, their {\em algebraic intersection number} is $0$.

The crossings between disjoint edges of $W$ and $W'$, counted with the signs of Section~\ref{sub:quads}, add up to $\Cr(C,C')$. The remaining contributions to their algebraic intersection number come from the common vertices of $W$ and $W'$, where the two curves meet. These contributions can be computed by ``pushing'' $W'$ slightly off each common vertex, and they only depend on the edges of $W$ and $W'$ at the vertex and on the rotation at that vertex. These intersections add up to $\Rot(C,C')$, and so the Cycle Identity $\Cr(C,C')=-\Rot(C,C')$ of Section~\ref{sec:cycproof} simply says that the total signed number of crossings is zero. This is essentially how the identity is proved for drawings in~\cite[Lemma~2.1]{openai}.


\subsection{The Cycle Identity}\label{sec:cycproof}

The Cycle Identity states that the quad count~\eqref{eq:J} and the rotation count~\eqref{eq:Rot} of any two closed walks, and more generally of any two balanced matrices, are negatives of each other.

\begin{theorem}[Cycle Identity]\label{thm:cancel}
For all balanced matrices $C,C'$,
$$
\Cr(C,C')=-\Rot(C,C').
$$
\end{theorem}

We will use the following standard linear algebra notion.

\begin{defbox}{Bilinear forms}
Let $V$ be a real vector space. A {\em bilinear form} on $V$ is a function $\beta$ that assigns a real number $\beta(a,b)$ to any two vectors $a,b\in V$, and that is linear in $a$ and in $b$. We say that $\beta$ is {\em antisymmetric} if $\beta(b,a)=-\beta(a,b)$ for all $a,b\in V$.
\end{defbox}

\vglue 0.3 cm
\noindent{\bf Remark.} {\em The word {\em antisymmetric} now has three meanings in this paper, all of the same kind: (i) an antisymmetric pair polynomial changes sign when $x_1$ and $x_2$ are exchanged, and also when $y_1$ and $y_2$ are exchanged; (ii) an antisymmetric matrix changes sign when rows and columns are exchanged; and now (iii) an antisymmetric form changes sign when its two arguments are exchanged. We note that this choice of terminology is standard in the literature.}
\vglue 0.3 cm

The proof of Theorem~\ref{thm:cancel} rests on one structural fact about antisymmetric bilinear forms on the space of balanced matrices (Lemma~\ref{lem:formzero} in Section~\ref{sub:structural}), and on a short computation of the rotation count at pairs of triangle matrices in a quad (Section~\ref{sub:cycleproof}). 

In the proof of the Cycle Identity we make essential use of the following two remarks, which follow easily from the definitions of $\Cr$ and $\Rot$ in~\eqref{eq:J} and~\eqref{eq:Rot}, respectively.

\begin{remark}\label{rem:bilinear}
Both $\Cr$ and $\Rot$ are bilinear forms on the space of antisymmetric $n\times n$ matrices.
\end{remark}

\begin{remark}\label{rem:antisym}
Both $\Cr$ and $\Rot$ are antisymmetric forms. That is, $\Cr(C',C)=$ $-\Cr(C,C')$ and $\Rot(C',C)=-\Rot(C,C')$ for all antisymmetric matrices $C$ and $C'$.
\end{remark}

\subsection{A structural lemma on  antisymmetric bilinear forms}\label{sub:structural}

For the rest of Section~\ref{sec:proofmainid} all the bilinear forms under consideration are antisymmetric. Moreover, for brevity, let us agree on the following.

\vglue 0.3 cm
\noindent{\bf Remark} (Forms). {\em For the rest of Section~\ref{sec:proofmainid} an antisymmetric bilinear form is simply called a {\em form}.}
\vglue 0.3 cm

We let $\mathcal B$ be the space of balanced $n\times n$ matrices.

The structural lemma of this section says that every form on $\mathcal B$ comes from a very simple kind of function, which we now introduce.

\begin{defbox}{Pair functions}
A {\em pair function} $\theta$ assigns a real number $\theta(ij,k\ell)$ to every ordered pair $(ij,k\ell)$ of directed edges with four distinct endvertices, in such a way that $\theta(k\ell,ij)=\theta(ji,k\ell)=\theta(ij,\ell k)=-\theta(ij,k\ell)$.

\medskip
Thus a pair function is determined by its three values $\theta(ij,k\ell)$, $\theta(ik,j\ell)$ and $\theta(i\ell,jk)$ on each quad $\{i,j,k,\ell\}$ with $i<j<k<\ell$, and these values can be chosen freely.
\end{defbox}

For instance, it is easy to verify that $\sign(ij,k\ell)$ is a pair function. As we will now see, each pair function gives a form, just as the pair function $\sign$ gives the quad count~\eqref{eq:J}.

\begin{defbox}{The form $\langle C,C'\rangle_\theta$ of a pair function $\theta$}
Let $\theta$ be any pair function. For antisymmetric matrices $C,C'$ we let
$$
\langle C,C'\rangle_\theta=\sum_{\substack{i<j,\ k<\ell\\ ij,\,k\ell\text{ disjoint edges}}}C(i,j)\,C'(k,\ell)\,\theta(ij,k\ell).
$$
\end{defbox}

In particular, $\langle C,C'\rangle_{\sign}$ is the quad count $\Cr(C,C')$ of~\eqref{eq:J}. Exactly as for $\Cr$, $\langle\cdot,\cdot\rangle_\theta$ is a form, and it depends linearly on $\theta$.

The next lemma is the key to the proof of the Cycle Identity. To motivate it we first recall that, as we have just noted, the quad count is the form of the pair function $\sign$. On the other hand, the rotation count is defined through pairs of edges with a common endvertex, and so there is no reason to expect it to be the form of any pair function. The lemma says that on balanced matrices every form is the form of exactly one pair function.  In particular, it follows that on balanced matrices the rotation count is also the form of a pair function. The Cycle Identity claims precisely that this pair function is $-\sign$. Since a pair function is determined by its values on each quad, this turns the Cycle Identity into a verification that it suffices to perform inside a single quad, a verification that we carry out in Section~\ref{sub:cycleproof}.
\begin{lemma}\label{lem:formzero}
Every form on $\mathcal B$ is the form $\langle\cdot,\cdot\rangle_\theta$ of exactly one pair function $\theta$.
\end{lemma}

\begin{proof}[Proof of Lemma~\ref{lem:formzero}]

Let $\mathcal M$ be the linear map from the space of pair functions to the space of forms on $\mathcal B$, defined by $\theta\mapsto \langle\cdot,\cdot\rangle_\theta$. We note that Lemma~\ref{lem:formzero} may be equivalently paraphrased by saying that {\em $\mathcal M$ is a bijection}.

Throughout the proof we say that two triangle matrices are {\em inside a quad} if the vertices of both triangles lie in a common quad. 

We prove two facts about $\mathcal M$:

\smallskip
\noindent(i)\enspace{\em If $\theta$ is a pair function such that $\langle T,T'\rangle_\theta=0$ for any two triangle matrices $T$ and $T'$ inside a quad, then $\theta$ is the zero pair function.}

\smallskip
\noindent(ii)\enspace{\em The space of pair functions and the space of forms on $\mathcal B$ have the same dimension, namely $3\binom{n}{4}$.}

\smallskip

We claim that the lemma (equivalently, that $\mathcal M$ is a bijection) follows from (i) and (ii). To see this we note that by~(i) the only pair function whose form is the zero form is the zero pair function (since the zero form is zero at every pair of triangle matrices inside a quad), and so $\mathcal M$ is injective. Now by~(ii) $\mathcal M$ is a linear map between two spaces of the same dimension, and so it is also surjective. Therefore (i) and (ii) indeed imply that $\mathcal M$ is a bijection, as claimed.

In view of this, we conclude the proof by showing (i) and (ii).

\medskip
\noindent{\em Proof of (i)}.\par\medskip

\noindent Suppose that $\langle T,T'\rangle_\theta=0$ for any two triangle matrices $T$ and $T'$ whose vertices lie in a common quad. The four endvertices of two disjoint edges form a quad, and so it suffices to show that $\theta$ is zero at every ordered pair of disjoint edges with endvertices in a quad $Q$.

Let $Q=\{i,j,k,\ell\}$ be a quad, with $i<j<k<\ell$, and let $\theta_1:=\theta(ij,k\ell)$, $\theta_2:=\theta(ik,j\ell)$ and $\theta_3:=\theta(i\ell,jk)$. As we noted in the box ``Pair functions'', these three numbers determine the value of $\theta$ at every ordered pair of disjoint edges with endvertices in $Q$. So our goal is to show that $\theta_1=\theta_2=\theta_3=0$.

Now $Q$ has four triangles, and so it has six pairs of triangle matrices. For each of these pairs the assumption is that the value of $\langle\cdot,\cdot\rangle_\theta$ at the pair is zero. This value is a combination of $\theta_1,\theta_2,\theta_3$, and so each pair gives a linear equation in $\theta_1,\theta_2,\theta_3$. We now exhibit three of these equations, and then show that their only solution is $\theta_1=\theta_2=\theta_3=0$. These equations come from the three triangles that contain $i$, namely $ijk$, $ij\ell$ and $ik\ell$, and they are shown in Table~\ref{tab:threepairs}.
\begin{table}[h!]
\centering
\renewcommand{\arraystretch}{1.25}
\begin{tabular}{c|c|c}
triangle matrices & pairs of their edges with four distinct endvertices & value of $\langle\cdot,\cdot\rangle_\theta$\\ \hline
$T_{ijk},\ T_{ij\ell}$ & $jk$ with $i\ell$, and $ik$ with $j\ell$ & $\theta_3-\theta_2$\\
$T_{ijk},\ T_{ik\ell}$ & $ij$ with $k\ell$, and $jk$ with $i\ell$ & $\theta_1+\theta_3$\\
$T_{ij\ell},\ T_{ik\ell}$ & $ij$ with $k\ell$, and $j\ell$ with $ik$ & $\theta_1-\theta_2$
\end{tabular}
\caption{The three pairs of triangle matrices of the quad $Q=\{i,j,k,\ell\}$ formed by the triangles that contain $i$, and the value of $\langle\cdot,\cdot\rangle_\theta$ at each of them (see~(i) in the proof of Lemma~\ref{lem:formzero}).}
\label{tab:threepairs}
\end{table}

We verify the first row, as the verification for the other two is totally~analogous. Above the diagonal $T_{ijk}$ has entries $1$ at $(i,j)$ and $(j,k)$, and $-1$ at $(i,k)$, and $T_{ij\ell}$ has entries $1$ at $(i,j)$ and $(j,\ell)$, and $-1$ at $(i,\ell)$. Thus $\langle T_{ijk},T_{ij\ell}\rangle_\theta=$ $T_{ijk}(j,k)\,T_{ij\ell}(i,\ell)\,\theta(jk,i\ell)$ $+T_{ijk}(i,k)\,T_{ij\ell}(j,\ell)\,\theta(ik,j\ell)=-\theta(jk,i\ell)-\theta(ik,j\ell)$, and by the sign rules of a pair function this is $\theta_3-\theta_2$, as claimed.

Now each row of Table~\ref{tab:threepairs} involves a pair of triangle matrices whose vertices lie in the quad $Q$. By our assumption that $\langle T,T'\rangle_\theta=0$ for any such pair, the three values in Table~\ref{tab:threepairs} are zero. That is, $\theta_3-\theta_2=0$, $\theta_1+\theta_3=0$ and $\theta_1-\theta_2=0$. The first and third equations imply that $\theta_1=\theta_2=\theta_3$, and then the second one gives $2\theta_1=0$. Therefore $\theta_1=\theta_2=\theta_3=0$. Since this holds for every quad we conclude that $\theta$ is the zero pair function, as claimed.

\medskip
\noindent{\em Proof of (ii)}.\par\medskip

\noindent We compute the two dimensions separately. The first one is straightforward. Indeed, a pair function is determined by three values on each quad, and so the pair functions form a space of dimension $3\binom{n}{4}$.

For the second one we first find the dimension of $\mathcal B$. For this we note that by Lemma~\ref{lem:span} the triangle matrices $T_{1jk}$ with $2\le j<k\le n$ span $\mathcal B$.

Moreover, we claim that they are linearly independent. Indeed, among the triangles $1j'k'$ with $2\le j'<k'\le n$, only $1jk$ contains the edge $jk$. Thus $T_{1jk}$ is the only one of these matrices whose entry in row $j$ and column $k$ is nonzero, and in a linear combination of them that equals the zero matrix the coefficient of $T_{1jk}$ must be zero.

Therefore these matrices form a basis of $\mathcal B$. Since there is one such matrix for each pair $j<k$ of elements of $\{2,\ldots,n\}$, we conclude that $\mathcal B$ has dimension $m=\binom{n-1}{2}$.

Now an antisymmetric bilinear form on a space of dimension $m$ is determined by its values at the pairs $(b_s,b_t)$ with $s<t$ of a basis $b_1,\ldots,b_m$, since by antisymmetry $(b_s,b_s)$ is mapped to $0$ and $(b_t,b_s)$ is mapped to the negative of the value at $(b_s,b_t)$. Since the values at the pairs $(b_s,b_t)$ with $s<t$ can be chosen arbitrarily, it follows that the forms on $\mathcal B$ form a space of dimension
$$
\binom{m}{2}=\frac{1}{2}\cdot\frac{(n-1)(n-2)}{2}\cdot\frac{n(n-3)}{2}=\frac{n(n-1)(n-2)(n-3)}{8}=3\binom{n}{4}.
$$
\end{proof}

\subsection{Proof of the Cycle Identity}\label{sub:cycleproof}

As we shall see in a moment, in order to prove the Cycle Identity it will suffice to show that $\Cr$ and $-\Rot$ agree at the three pairs of triangle matrices of Table~\ref{tab:threepairs}, for every quad. With this in mind we start by computing $\Rot$ at these pairs.

\begin{lemma}[The rotation count at the triangle matrices of a quad]\label{lem:rotquad}
Let $Q=\{i,j,k,\ell\}$ be a quad, with $i<j<k<\ell$, and let $\orient{i},\orient{j},\orient{k},\orient{\ell}$ be the orientations of its four vertices. Then
$$
\Rot(T_{ijk},T_{ij\ell})=\tfrac12\,(\orient{i}+\orient{j}),\qquad
\Rot(T_{ijk},T_{ik\ell})=\tfrac12\,(\orient{i}-\orient{k}),\qquad
\Rot(T_{ij\ell},T_{ik\ell})=\tfrac12\,(\orient{i}+\orient{\ell}).
$$
\end{lemma}

\begin{proof}[Proof of Lemma~\ref{lem:rotquad}]
The three triangle matrices of the lemma are shown in Figure~\ref{fig:trimatrices}.

\begin{figure}[ht!]
\centering
$
\begin{array}{c|cccc}
T_{ijk} & \hphantom{-}i & \hphantom{-}j & \hphantom{-}k & \hphantom{-}\ell\\ \hline
i & \hphantom{-}0 & \hphantom{-}1 & -1 & \hphantom{-}0\\
j & -1 & \hphantom{-}0 & \hphantom{-}1 & \hphantom{-}0\\
k & \hphantom{-}1 & -1 & \hphantom{-}0 & \hphantom{-}0\\
\ell & \hphantom{-}0 & \hphantom{-}0 & \hphantom{-}0 & \hphantom{-}0
\end{array}
\qquad\quad
\begin{array}{c|cccc}
T_{ij\ell} & \hphantom{-}i & \hphantom{-}j & \hphantom{-}k & \hphantom{-}\ell\\ \hline
i & \hphantom{-}0 & \hphantom{-}1 & \hphantom{-}0 & -1\\
j & -1 & \hphantom{-}0 & \hphantom{-}0 & \hphantom{-}1\\
k & \hphantom{-}0 & \hphantom{-}0 & \hphantom{-}0 & \hphantom{-}0\\
\ell & \hphantom{-}1 & -1 & \hphantom{-}0 & \hphantom{-}0
\end{array}
\qquad\quad
\begin{array}{c|cccc}
T_{ik\ell} & \hphantom{-}i & \hphantom{-}j & \hphantom{-}k & \hphantom{-}\ell\\ \hline
i & \hphantom{-}0 & \hphantom{-}0 & \hphantom{-}1 & -1\\
j & \hphantom{-}0 & \hphantom{-}0 & \hphantom{-}0 & \hphantom{-}0\\
k & -1 & \hphantom{-}0 & \hphantom{-}0 & \hphantom{-}1\\
\ell & \hphantom{-}1 & \hphantom{-}0 & -1 & \hphantom{-}0
\end{array}
$
\caption{The triangle matrices $T_{ijk}$, $T_{ij\ell}$ and $T_{ik\ell}$ of a quad $Q=\{i,j,k,\ell\}$. We only show the rows and columns $i,j,k,\ell$, since all the other entries are zero.}
\label{fig:trimatrices}
\end{figure}

We note that each of the equations involves two triangles that share two vertices. Now a glance at the definition of $\Rot$ (see~\eqref{eq:Rot}) shows that only the two common vertices of the two triangles contribute to each rotation count. Using that $\tau_m(j,j)=0$ and $\tau_m(j,k)=-\tau_m(k,j)$, a straightforward calculation yields
\begin{align*}
\Rot(T_{ijk},T_{ij\ell})&{}=\tfrac12\big(\tau_i(j,k)+\tau_i(k,\ell)+\tau_i(\ell,j)\big)+\tfrac12\big(\tau_j(i,k)+\tau_j(k,\ell)+\tau_j(\ell,i)\big),\\
\Rot(T_{ijk},T_{ik\ell})&{}=\tfrac12\big(\tau_i(j,k)+\tau_i(k,\ell)+\tau_i(\ell,j)\big)-\tfrac12\big(\tau_k(i,j)+\tau_k(j,\ell)+\tau_k(\ell,i)\big),\\
\Rot(T_{ij\ell},T_{ik\ell})&{}=\tfrac12\big(\tau_i(j,k)+\tau_i(k,\ell)+\tau_i(\ell,j)\big)+\tfrac12\big(\tau_\ell(i,j)+\tau_\ell(j,k)+\tau_\ell(k,i)\big).
\end{align*}
For instance, row $i$ of $T_{ijk}$ has entry $1$ at $j$ and $-1$ at $k$, and row $i$ of $T_{ij\ell}$ has entry $1$ at $j$ and $-1$ at $\ell$. Thus by~\eqref{eq:Rot} the contribution of $i$ to $\Rot(T_{ijk},T_{ij\ell})$ is
$$
\tfrac12\big(\tau_i(j,j)-\tau_i(j,\ell)-\tau_i(k,j)+\tau_i(k,\ell)\big)=\tfrac12\big(\tau_i(j,k)+\tau_i(k,\ell)+\tau_i(\ell,j)\big),
$$
and the contribution of $j$ to $\Rot(T_{ijk},T_{ij\ell})$ is calculated similarly. The equalities for $\Rot(T_{ijk},T_{ik\ell})$ and $\Rot(T_{ij\ell},T_{ik\ell})$ are handled in a totally analogous manner.

In order to finish the proof it clearly suffices to show that
\begin{align*}
\tau_i(j,k)+\tau_i(k,\ell)+\tau_i(\ell,j)&{}=\orient{i}, & \tau_j(i,k)+\tau_j(k,\ell)+\tau_j(\ell,i)&{}=\orient{j},\\
\tau_k(i,j)+\tau_k(j,\ell)+\tau_k(\ell,i)&{}=\orient{k}, & \tau_\ell(i,j)+\tau_\ell(j,k)+\tau_\ell(k,i)&{}=\orient{\ell}.
\end{align*}

\begin{figure}[ht!]
\centering
$
\begin{array}{c|ccc|c|c}
\text{order in }\pi(i) & \tau_i(j,k) & \tau_i(k,\ell) & \tau_i(\ell,j) & \text{sum} & \orient{i}\\ \hline
j,k,\ell & +1 & +1 & -1 & +1 & +1\\
k,\ell,j & -1 & +1 & +1 & +1 & +1\\
\ell,j,k & +1 & -1 & +1 & +1 & +1\\
\hline
\ell,k,j & -1 & -1 & +1 & -1 & -1\\
k,j,\ell & -1 & +1 & -1 & -1 & -1\\
j,\ell,k & +1 & -1 & -1 & -1 & -1\\
\end{array}
$
\caption{The six orders in which $j<k<\ell$ can appear in $\pi(i)$. In the first three the rotation at $i$ in $\Pi|_Q$ is increasing, and in the last three it is decreasing. In every row the sum of the three rotation signs is $\orient{i}$.}
\label{fig:sixorders}
\end{figure}

We work out the calculation only for the first equality, as the other three are handled in a totally analogous manner. The vertices $j,k,\ell$ appear in $\pi(i)$ in one of six orders, which are listed in Figure~\ref{fig:sixorders}. For each of the six orders we show in the corresponding row the three terms of the sum. In every case the sum agrees with $\orient{i}$, as claimed.
\end{proof}

\begin{proof}[Proof of Theorem~\ref{thm:cancel}]
Let $Q=\{i,j,k,\ell\}$ be any quad, with $i<j<k<\ell$. We claim that in order to prove the theorem it suffices to show that
\setlength{\abovedisplayskip}{10pt}\setlength{\belowdisplayskip}{10pt}\setlength{\abovedisplayshortskip}{10pt}\setlength{\belowdisplayshortskip}{10pt}
\begin{equation*}\tag{$*$}
\Cr(C,C')+\Rot(C,C')=0\quad\text{for}\quad (C,C')=(T_{ijk},T_{ij\ell}),\ (T_{ijk},T_{ik\ell})\ \text{and}\ (T_{ij\ell},T_{ik\ell}).
\end{equation*}
These are the three pairs of triangle matrices displayed in Table~\ref{tab:threepairs}.

To see why it suffices to show ($*$), we first recall that Remarks~\ref{rem:bilinear} and~\ref{rem:antisym} imply that $\Cr+\Rot$ is a form on $\mathcal B$, and by Lemma~\ref{lem:formzero} it follows that $\Cr+\Rot$ must be $\langle\cdot,\cdot\rangle_\theta$ for some pair function $\theta$. Assuming that ($*$) holds we have
\begin{equation}\label{eq:threezero}
\langle T_{ijk},T_{ij\ell}\rangle_\theta=\langle T_{ijk},T_{ik\ell}\rangle_\theta=\langle T_{ij\ell},T_{ik\ell}\rangle_\theta=0.
\end{equation}

We recall that in Table~\ref{tab:threepairs} we write $\theta_1,\theta_2$, and $\theta_3$ for $\theta(ij,k\ell)$, $\theta(ik,j\ell)$, and $\theta(i\ell,jk)$, respectively. Thus Table~\ref{tab:threepairs} yields
\begin{equation}\label{eq:threetheta}
\begin{aligned}
\langle T_{ijk},T_{ij\ell}\rangle_\theta&{}=\theta(i\ell,jk)-\theta(ik,j\ell),\\
\langle T_{ijk},T_{ik\ell}\rangle_\theta&{}=\theta(ij,k\ell)+\theta(i\ell,jk), \quad\text{and} \\
\langle T_{ij\ell},T_{ik\ell}\rangle_\theta&{}=\theta(ij,k\ell)-\theta(ik,j\ell).
\end{aligned}
\end{equation}
Using~\eqref{eq:threezero} and~\eqref{eq:threetheta} we obtain $\theta(ij,k\ell)=\theta(ik,j\ell)=\theta(i\ell,jk)=0$. Since $Q$ is arbitrary, it follows that $\theta$ is the zero pair function, and so $\Cr=-\Rot$, as claimed by the theorem. Thus in order to prove the theorem it suffices to show ($*$).

To show ($*$) we recall from Section~\ref{sub:structural} that $\Cr=\langle\cdot,\cdot\rangle_{\sign}$. We claim that
\begin{align*}
\Cr(T_{ijk},T_{ij\ell})&{}=\sign(i\ell,jk)-\sign(ik,j\ell)=-\tfrac12(\orient{i}+\orient{j})=-\Rot(T_{ijk},T_{ij\ell}),\\
\Cr(T_{ijk},T_{ik\ell})&{}=\sign(ij,k\ell)+\sign(i\ell,jk)=-\tfrac12(\orient{i}-\orient{k})=-\Rot(T_{ijk},T_{ik\ell}),\\
\Cr(T_{ij\ell},T_{ik\ell})&{}=\sign(ij,k\ell)-\sign(ik,j\ell)=-\tfrac12(\orient{i}+\orient{\ell})=-\Rot(T_{ij\ell},T_{ik\ell}).
\end{align*}
Indeed, in each line the first equality follows from Table~\ref{tab:threepairs} with $\theta=\sign$, since $\Cr=\langle\cdot,\cdot\rangle_{\sign}$, the second equality follows from~\eqref{eq:signs}, and the third follows from Lemma~\ref{lem:rotquad}. These three equalities prove~($*$), and so we are done.
\end{proof}


\section{Proof of the Main Identity II: adding the polynomials}\label{sec:polysec}
\enlargethispage{\baselineskip}

We now show that the Main Identity (Theorem~\ref{thm:master}) follows easily from the Cycle Identity (Theorem~\ref{thm:cancel}). For the reader's convenience, we restate the theorem.

\theoremstyle{plain}
\newtheorem*{restatedmaster}{Theorem~\ref{thm:master}}
\begin{restatedmaster}[Main Identity]
For every antisymmetric pair polynomial $A$,
\begin{equation*}
\sum_{(ij,k\ell)\in\xx}\sign(ij,k\ell)\,\omega_i\omega_j\omega_k\omega_\ell\,A(i,j,k,\ell)
\;=\;
-\tfrac{1}{2}\sum_{i}\omega_i^2\sum_{j,k\ne i}\tau_i(j,k)\,\omega_j\omega_k\,A(i,j,i,k).
\tag{\ref{eq:master}}
\end{equation*}
\end{restatedmaster}

In the proof of the theorem and in a later section it will be very convenient to use the following notion.

\begin{definition}[Antisymmetric version of a polynomial]\label{def:antiversion}
Let $H=H(x_1,x_2,y_1,y_2)$ be a pair polynomial. The {\em antisymmetric version} of $H$ is the polynomial ${\widehat H}$ defined by
\begin{equation*}
{\widehat H}(x_1,x_2,y_1,y_2)=\tfrac{1}{4}\big(H(x_1,x_2,y_1,y_2)-H(x_2,x_1,y_1,y_2)-H(x_1,x_2,y_2,y_1)+H(x_2,x_1,y_2,y_1)\big).
\end{equation*}
\end{definition}

It is easy to verify that ${\widehat H}$ is always an antisymmetric pair polynomial.

In a nutshell, the proof of the Main Identity consists of two observations. First, every antisymmetric pair polynomial is a linear combination of the antisymmetric versions of monomials. Since both sides of~\eqref{eq:master} are linear in $A$, it follows that it suffices to prove the Main Identity for the antisymmetric versions of monomials. Second, for the antisymmetric version of a monomial the Main Identity is exactly the Cycle Identity for two particular balanced matrices $C$ and $C'$.

The proof relies on the next lemma, which involves the weights $\omega_i$ defined in Section~\ref{sec:conv}. For easy reference, we remind the reader that
$$
\omega_i=\frac{1}{\prod_{j\in[n]\setminus\{i\}}(i-j)}.
$$

\begin{lemma}\label{lem:lagrange}
For every polynomial $G$ of degree at most $n-2$,
$$
\sum_{i=1}^n\omega_i\,G(i)=0.
$$
\end{lemma}

\begin{proof}[Proof of Lemma~\ref{lem:lagrange}]
Let $G$ be a polynomial of degree at most $n-2$. By Lagrange interpolation (see for instance \cite[Theorem~3.2]{burden-faires}) $G$ satisfies
$$
G(t)=\sum_{i=1}^nG(i)\prod_{j\ne i}\frac{t-j}{i-j}=\sum_{i=1}^n\omega_i\,G(i)\prod_{j\ne i}(t-j),
$$
since both sides have degree at most $n-1$ and agree at the $n$ points $1,\ldots,n$. Now each product $\prod_{j\ne i}(t-j)$ has $n-1$ factors of the form $t-j$, and so its leading term is $t^{n-1}$. The coefficient of $t^{n-1}$ on the right-hand side is $\sum_i\omega_i\,G(i)$. Finally, since the left-hand side is $G$ and $\deg G\le n-2$, it follows that this coefficient is $0$.
\end{proof}

With Lemma~\ref{lem:lagrange} at hand we can introduce the balanced matrices used in the proof. For any two polynomials $F_1$ and $F_2$ in one variable, each of them of degree at most $n-2$, we let $C_{F_1,F_2}$ be the $n\times n$ matrix with entries
\begin{equation}\label{eq:Cfg}
C_{F_1,F_2}(i,j)=\tfrac{1}{2}\,\omega_i\omega_j\big(F_1(i)F_2(j)-F_1(j)F_2(i)\big).
\end{equation}

We note that for any two polynomials $F_1$ and $F_2$ the matrix $C_{F_1,F_2}$ is clearly antisymmetric. Moreover, it is also balanced, since the entries of row $i$ add up to
$$
\tfrac{1}{2}\,\omega_i\Big(F_1(i)\sum_{j}\omega_j\,F_2(j)-F_2(i)\sum_{j}\omega_j\,F_1(j)\Big),
$$
and Lemma~\ref{lem:lagrange} implies that both sums are $0$.

\enlargethispage{\baselineskip}
\begin{proof}[Proof of Theorem~\ref{thm:master}]
The proof has two steps. We show that (i) Equation~\eqref{eq:master} holds for $A=\widehat H$, where $\widehat H$ is the antisymmetric version of any monomial $H=x_1^{s_1}x_2^{s_2}y_1^{t_1}y_2^{t_2}$, where $0 \le s_1, s_2, t_1, t_2 \le n-2$, and then we show that (ii) every antisymmetric pair polynomial is a linear combination of these polynomials $\widehat H$. Since both sides of~\eqref{eq:master} are linear in $A$, it clearly follows that (i) and (ii) imply the theorem.

To prove (i) we start by letting $F_1(x)=x^{s_1}$, $F_2(x)=x^{s_2}$, $F_3(x)=x^{t_1}$ and $F_4(x)=x^{t_2}$, so that $H(x_1,x_2,y_1,y_2)=F_1(x_1)\,F_2(x_2)\,F_3(y_1)\,F_4(y_2)$. We claim that~\eqref{eq:master} for $A={\widehat H}$ is exactly the Cycle Identity for the balanced matrices $C_{F_1,F_2}$ and $C_{F_3,F_4}$.

Indeed, a straightforward calculation using~\eqref{eq:Cfg} shows that $C_{F_1,F_2}(i,j)$ $\,C_{F_3,F_4}(k,\ell)=$ $\omega_i\omega_j\omega_k\omega_\ell\,{\widehat H}(i,j,k,\ell)$ for all $i,j,k,\ell$. Therefore $\Cr(C_{F_1,F_2},C_{F_3,F_4})$ is precisely the left-hand side of~\eqref{eq:master} for $A={\widehat H}$, and $-\Rot(C_{F_1,F_2},C_{F_3,F_4})$ is precisely its right-hand side. Now since the matrices $C_{F_1,F_2}$ and $C_{F_3,F_4}$ are balanced (as we remarked right after~\eqref{eq:Cfg}) the Cycle Identity (Theorem~\ref{thm:cancel}) applies to them, and therefore~\eqref{eq:master} holds for $A={\widehat H}$.

To prove (ii) let $A$ be any antisymmetric pair polynomial. Since $A$ is a pair polynomial it can be written as
$$
A=\sum_H c_H\,H,
$$
where the summation is over all monomials $H=x_1^{s_1}x_2^{s_2}y_1^{t_1}y_2^{t_2}$ with $0\le s_1,s_2,t_1,t_2\le n-2$, and the $c_H$ are real numbers.

Taking the antisymmetric version is a linear operation, and so
$$
{\widehat A}=\sum_H c_H\,{\widehat H}.
$$

Finally, since $A$ is antisymmetric, each of the four terms in Definition~\ref{def:antiversion} is equal to $A$, and so ${\widehat A}=A$. Therefore $A=\sum_H c_H\,{\widehat H}$, and so (ii) follows.
\end{proof}


\section{Proof of Proposition~\ref{pro:inj}}\label{app:inj}

We recall that Proposition~\ref{pro:inj} states that if $P(x_1,x_2,y_1,y_2)\in\ss$ is a polynomial that satisfies $P(i,j,k,\ell)=0$ for every $(ij,k\ell)\in\xx$, then $P$ is the zero polynomial. We recall that $P\in\ss$ means that $P$ has degree at most $r-2$ in each variable, and its value does not change when $x_1$ and $x_2$ are exchanged and also when $y_1$ and $y_2$ are exchanged.

The proof relies on Lemma~\ref{lem:detector} below. This lemma is a technical fact about polynomials, which involves only the weights $\omega_j$ and a ``sign'' function $\rho$ on pairs of vertices. Its proof involves several facts and arguments on polynomial algebra, and we defer it to Section~\ref{sec:interp}.

The sign function in the lemma comes from a linear order of the vertices other than a fixed vertex $i$. Given a linear order on $[n]\setminus\{i\}$, its {\em sign function} $\rho$ is defined by $\rho(j,k)=+1$ if $j$ occurs after $k$, $\rho(j,k)=-1$ if $j$ occurs before $k$, and $\rho(j,j)=0$. 

We note that the rotation sign $\tau_i$ is precisely a sign function of this kind. Indeed, $\tau_i$ is the sign function of the reverse of the linear order $\pi(i)$, and we will apply the lemma to $\rho=\tau_i$.

\begin{lemma}\label{lem:detector}
Let $i\in[n]$, let $\rho$ be the sign function of a linear order of the $2r$ integers in $[n]\setminus\{i\}$, and let $F(x_2,y_2)$ be a polynomial in two variables of degree at most $r-2$ in each variable. If
\begin{equation}\label{eq:detector}
\sum_{j,k\in[n]\setminus\{i\}}\rho(j,k)\,\omega_j\omega_k\,(j-k)(j-i)(k-i)\,F(j,k)\,B(j,k)=0
\end{equation}
for every polynomial $B(x_2,y_2)$ of degree at most $r-1$ in each variable, then $F$ is the zero polynomial.
\end{lemma}

We note that Lemma~\ref{lem:detector} is a  statement about integers and polynomials. The integers in $[n]$ naturally correspond to the vertices of $K_n$, but the lemma does not involve the rotation system $\Pi$.

\begin{proof}[Proof of Proposition~\ref{pro:inj}]
We proceed by contradiction. Suppose that there is a nonzero polynomial $P\in\ss$ such that $P(i,j,k,\ell)=0$ for every $(ij,k\ell)\in\xx$. Among all such polynomials we choose $P$ so that its degree in $x_1$ is as small as possible. Our goal is to arrive at a contradiction, by exhibiting another nonzero polynomial in $\ss$ that also vanishes at every element of $\xx$, and whose degree in $x_1$ is strictly smaller than the degree of $P$ in $x_1$ (this is Step 4 in the proof).

This contradiction proves the proposition, since it shows that the only polynomial in $\ss$ that vanishes at all the elements of $\xx$ is the zero polynomial, that is, that $E$ is injective.

We divide the proof into four steps. For convenience and easy reference we define right away some key objects that are used throughout the proof (see bulleted list below).

These objects have a clear purpose. In Step 1 we apply the Main Identity to $A_i=P\cdot R_i$. Since $P$ vanishes at every element of $\xx$ it follows that the quads side is zero, and so by the Main Identity the rotations side is zero too. {\em The objects are designed so that this rotations side involves only the values $P(i,j,i,k)$ at one prescribed vertex $i$, in exactly the form required by Lemma~\ref{lem:detector}.} We discuss these objects in more detail immediately after~\eqref{eq:Si}.

\vspace{6pt}
{\setlength{\leftmargini}{1.5em}\begin{itemize}\setlength{\itemsep}{8pt}
\item We let $B=B(x_2,y_2)$ be an arbitrary, unspecified polynomial in two variables, of degree at most $r-1$ in each variable.

\medskip
\item For each vertex $i$ we let $U_i$ and $V_i$ be polynomials in one variable, each of degree at most $r$, such that $U_i(i)V_i(i)=1$ and $U_i(j)V_i(j)=0$ for every vertex $j\ne i$.

\medskip
To see that $U_i$ and $V_i$ with these properties exist let $j_1<j_2<\cdots<j_{2r}$ be the $2r=n-1$ vertices other than $i$, and let
$$
U_i(t)=\omega_i\,(t-j_1)(t-j_2)\cdots(t-j_r)\quad\text{and}\quad V_i(t)=(t-j_{r+1})(t-j_{r+2})\cdots(t-j_{2r}).
$$
It is easy to see that $U_i$ and $V_i$ satisfy the required properties. Indeed, every vertex $j\ne i$ is a root of $U_i$ or of $V_i$, and $U_i(i)V_i(i)=\omega_i\prod_{j\ne i}(i-j)=1$ by the definition of $\omega_i$. The reader may wonder why we did not simply define $U_i$ and $V_i$ by these formulas, and the reason is that we prefer the bulleted description because it is lighter and states exactly the two properties that we use.

\medskip
\item For each vertex $i$, we let
\begin{equation}\label{eq:Hi}
H_i=(x_2-y_2)(x_2-y_1)(y_2-x_1)\cdot U_i(x_1)\,V_i(y_1)\cdot B(x_2,y_2),
\end{equation}
and we let $R_i$ be the antisymmetric version of $H_i$ (Definition~\ref{def:antiversion}), that is,
\begin{equation}\label{eq:Ri}
\begin{split}
R_i(x_1,x_2,y_1,y_2)\mathrel{:=}\tfrac{1}{4}\big(&H_i(x_1,x_2,y_1,y_2)-H_i(x_2,x_1,y_1,y_2)\\
&-H_i(x_1,x_2,y_2,y_1)+H_i(x_2,x_1,y_2,y_1)\big).
\end{split}
\end{equation}

\item For each vertex $i$, we let $\sigma_i$ be the number
\begin{equation}\label{eq:Si}
\sigma_i\mathrel{:=}\sum_{j,k\ne i}\tau_i(j,k)\,\omega_j\omega_k\,(j-k)(j-i)(k-i)\,P(i,j,i,k)\,B(j,k).
\end{equation}
\end{itemize}}
\bigskip

It seems worth explaining a bit more about the roles of these objects. As we mentioned above, we will apply the Main Identity to $A_i=P\cdot R_i$. We need $R_i$ to be antisymmetric in order to make $A_i$ antisymmetric, because the Main Identity applies to antisymmetric polynomials. We force the antisymmetry of $R_i$ in its definition~\eqref{eq:Ri}, and at the end of the day $R_i$ inherits the properties of $H_i$, which is a polynomial carefully crafted, as we now explain.

Each part of $H_i$ has one purpose. On the rotations side the variables $x_1$ and $y_1$ take the same value. There the three factors $(x_2-y_2)(x_2-y_1)(y_2-x_1)$ make three of the four terms of $R_i$ vanish, and only the first term (which is $H_i$ itself) remains. Now to see the roles of $U_i$ and $V_i$ we note that, by their very definitions, in $H_i$ the product $U_i(x_1)V_i(y_1)$ is zero when $x_1=y_1$ is a vertex other than $i$, and so only vertex $i$ survives on the rotations side. Finally $B$ is free, matching the freedom explicitly stated in Lemma~\ref{lem:detector}.

We finally note that the number $\sigma_i$ is simply the term we use for the sum that the rotations side becomes when we apply the Main Identity to $A_i=P\cdot R_i$.
\bigskip

Having defined these key objects we proceed to the rest of the proof, broken into four steps to help comprehension.

\bigskip

\noindent{\sc Step 1}.\enspace{\em For every vertex $i$, the number $\sigma_i$ defined in~\eqref{eq:Si} is zero.}\par\medskip

Let $i$ be a fixed vertex, and let $A_i=P\cdot R_i$. We claim that the Main Identity can be applied to $A_i$, that is, that $A_i$ is an antisymmetric pair polynomial.

To see that $A_i$ is antisymmetric we start by noting that $R_i$ is antisymmetric, since it is the antisymmetric version of the polynomial $H_i$ (Definition~\ref{def:antiversion}). Since the product of a symmetric polynomial (namely $P$) and an antisymmetric polynomial (namely $R_i$) is antisymmetric, it follows that $A_i$ is antisymmetric.

Now to see that $A_i$ is a pair polynomial, that is, it has degree at most $n-2$ in each variable, we first note that it is easy to verify that $H_i$ has degree at most $r+1$ in each variable, and from this it follows immediately that $R_i$ also has degree at most $r+1$ in each variable. Since $P$ has degree at most $r-2$ in each variable, it follows that $A_i$ has degree at most $(r-2)+(r+1)=2r-1=n-2$ in each variable. Thus $A_i$ is a pair polynomial, as claimed.

We now note that the left-hand side (the quads side) of the Main Identity applied to $A_i$ is zero. This simply follows because each term in the summation involves a value $A_i(j,k,\ell,m)=P(j,k,\ell,m)\,R_i(j,k,\ell,m)$ at an element $(jk,\ell m)$ of $\xx$, and $P$ vanishes at every element of $\xx$. 

Thus the quads side is zero, and so Theorem~\ref{thm:master} then implies that the rotations side is also zero, that is,
\begin{equation}\label{eq:rotzeroF}
\tfrac{1}{2}\sum_m\omega_m^2\sum_{j,k\ne m}\tau_m(j,k)\,\omega_j\omega_k\,P(m,j,m,k)\,R_i(m,j,m,k)=0.
\end{equation}

We now simplify the terms $R_i(m,j,m,k)$ in this equation. By the definition of $R_i$ in~\eqref{eq:Ri} we have
$$
4R_i(m,j,m,k)=H_i(m,j,m,k)-H_i(j,m,m,k)-H_i(m,j,k,m)+H_i(j,m,k,m).
$$
The last three terms are zero, since $H_i$ has the factor $x_2-y_1$, which vanishes at $(j,m,m,k)$, the factor $y_2-x_1$, which vanishes at $(m,j,k,m)$, and the factor $x_2-y_2$, which vanishes at $(j,m,k,m)$. Thus
\begin{equation}\label{eq:Rivalues}
R_i(m,j,m,k)=\tfrac{1}{4}H_i(m,j,m,k)=\tfrac{1}{4}(j-k)(j-m)(k-m)\,U_i(m)V_i(m)\,B(j,k).
\end{equation}
Now recall that $U_i(m)V_i(m)=0$ for every vertex $m\ne i$, and that $U_i(i)V_i(i)=1$. Thus $R_i(m,j,m,k)=0$ unless $m=i$, and so the only terms that are not zero in~\eqref{eq:rotzeroF} are those with $m=i$. Therefore~\eqref{eq:rotzeroF} implies that
\begin{equation}\label{eq:rotzeroi}
\tfrac{1}{2}\,\omega_i^2\sum_{j,k\ne i}\tau_i(j,k)\,\omega_j\omega_k\,P(i,j,i,k)\,R_i(i,j,i,k)=0.
\end{equation}
Using~\eqref{eq:Rivalues} with $m=i$, and also using that $U_i(i)V_i(i)=1$, it follows that $R_i(i,j,i,k)=\tfrac{1}{4}(j-k)(j-i)(k-i)\,B(j,k)$. Substituting this into~\eqref{eq:rotzeroi} we obtain
$$
\tfrac{1}{8}\,\omega_i^2\sum_{j,k\ne i}\tau_i(j,k)\,\omega_j\omega_k\,(j-k)(j-i)(k-i)\,P(i,j,i,k)\,B(j,k)=0.
$$
The sum in this equation is exactly the number $\sigma_i$ defined in~\eqref{eq:Si}. Thus $\tfrac{1}{8}\omega_i^2\,\sigma_i=0$. Since $\omega_i\ne 0$ for every vertex $i$ it follows that $\sigma_i=0$.

\medskip
\noindent{\sc Step 2}.\enspace{\em For every vertex $i$, the polynomial $P(i,x_2,i,y_2)$ is the zero polynomial.}\par\medskip

Let $i$ be any fixed vertex, and let $F(x_2,y_2)\mathrel{:=}P(i,x_2,i,y_2)$. We claim that the hypotheses of Lemma~\ref{lem:detector} hold for this $F$, with $\rho=\tau_i$.

To see this first we note that $\tau_i$ is indeed the sign function of a linear order of the vertices other than $i$, namely the reverse of the rotation $\pi(i)$. Next we note that since $P$ has degree at most $r-2$ in each variable then so does $F$. Finally we note that with $\rho=\tau_i$ and $F=P(i,x_2,i,y_2)$, the sum in~\eqref{eq:detector} is the number $\sigma_i$ defined in~\eqref{eq:Si}. Since $B$ is an arbitrary polynomial of degree at most $r-1$ in each variable, Step~1 shows that this sum is zero for every such $B$, that is, it shows that~\eqref{eq:detector} holds for every such $B$. Thus the hypotheses of Lemma~\ref{lem:detector} indeed hold for $F$, with $\rho=\tau_i$, as claimed.

Applying Lemma~\ref{lem:detector} we then obtain that $F$ is the zero polynomial, as claimed.

\medskip
\noindent{\sc Step 3}.\enspace{\em There is a polynomial ${\widetilde P}$ such that $P=M\cdot{\widetilde P}$, where $M=(x_1-y_1)(x_1-y_2)(x_2-y_1)(x_2-y_2)$.}\par\medskip

We start by considering the polynomial in three variables
$$
G(x_1,x_2,y_2)\mathrel{:=}P(x_1,x_2,x_1,y_2),
$$
which is obtained from $P=P(x_1,x_2,y_1,y_2)$ by making $y_1=x_1$. By Step 2 the polynomial $G(i,x_2,y_2)$ is the zero polynomial for every vertex $i$. Now $G$ has degree at most $2(r-2)=2r-4<n$ in $x_1$, and so for any fixed values of $x_2$ and $y_2$ we may regard $G$ as a polynomial $J(x_1)$ in $x_1$ of degree less than $n$, and since $1,\ldots,n$ are all distinct roots of $J(x_1)$ it follows that $J$ must be the zero polynomial. 

Therefore $G(x_1,x_2,y_2)=0$ for all values of $x_1$, $x_2$ and $y_2$, that is, $G$ also is the zero polynomial. We claim that this implies that $x_1-y_1$ divides $P$. To see this we note that if we divide $P$ by $x_1-y_1$, regarding both as polynomials in $y_1$, we obtain $P=(x_1-y_1)\,P_1+P_0$ for some polynomials $P_1$ and $P_0$, where $P_0$ does not involve $y_1$. Setting $y_1=x_1$ gives $P_0=G$, and so $P_0$ is the zero polynomial. Thus $P=(x_1-y_1)\,P_1$, and so $x_1-y_1$ divides $P$, as claimed.

Since $P$ does not change when $x_1$ and $x_2$ are exchanged, and it does not change either when $y_1$ and $y_2$ are exchanged, it follows that $x_2-y_1$, $x_1-y_2$ and $x_2-y_2$ must also divide $P$. These four polynomials are irreducible, and no two of them are constant multiples of each other, and so by unique factorization their product $M$ divides $P$. In other words, there is a polynomial ${\widetilde P}$ such that $P=M\cdot{\widetilde P}$.

\medskip
\noindent{\sc Step 4}.\enspace{\em The polynomial ${\widetilde P}$ from Step 3 is nonzero, belongs to $\ss$, vanishes at every element of $\xx$, and has smaller degree in $x_1$ than $P$. This contradicts the choice of $P$ at the beginning of the proof of the proposition.}\par\medskip

We start by noting that ${\widetilde P}$ is nonzero, because $P=M\cdot {\widetilde P}$ and by assumption $P$ is not the zero polynomial.

Now to see that ${\widetilde P}$ belongs to $\ss$ we first recall that this means that ${\widetilde P}$ does not change when $x_1$ and $x_2$ are exchanged, does not change when $y_1$ and $y_2$ are exchanged, and that it has degree at most $r-2$ in each variable. Exchanging $x_1$ and $x_2$ leaves $M$ unchanged, and the same holds for $y_1$ and $y_2$. Since $P=M\cdot{\widetilde P}$ has these two symmetries then ${\widetilde P}$ also has these symmetries. Finally, since $M$ has degree exactly $2$ in each variable it follows that the degree of ${\widetilde P}$ in each variable is two less than the degree of $P$ in that variable, and so ${\widetilde P}$ has degree at most $r-2$ in each variable. Thus ${\widetilde P}\in\ss$, as claimed. 

We now note that Step 3 implies that the degree of ${\widetilde P}$ in $x_1$ is two less than the degree of $P$ in $x_1$, and so it only remains to show that ${\widetilde P}$ vanishes at every element of $\xx$. To see this let $(ij,k\ell)$ be any element of $\xx$. The four vertices $i,j,k,\ell$ are distinct, and so $M(i,j,k,\ell)\ne 0$. Now since $P(i,j,k,\ell)=M(i,j,k,\ell)\,{\widetilde P}(i,j,k,\ell)=0$ it follows that ${\widetilde P}(i,j,k,\ell)=0$, and so we are done.

Thus ${\widetilde P}$ satisfies the properties in Step 4. Since its degree in $x_1$ is smaller than the degree of $P$ in $x_1$ this contradicts the choice of $P$, and so the proof of the proposition is complete.
\end{proof}

\section{Proof of Proposition~\ref{pro:iso}}\label{sec:iso}

We recall that Proposition~\ref{pro:iso} states that $\dim E(\ss)\le\tfrac{1}{2}|\xx|$. For easy reference, let us recall the three objects involved in this inequality.

First, $\ss$ is the space of polynomials $P(x_1,x_2,y_1,y_2)$ of degree at most $r-2$ in each variable that do not change when $x_1$ and $x_2$ are exchanged, nor when $y_1$ and $y_2$ are exchanged. Second, $\xx$ is the set of ordered pairs $(ij,k\ell)$ of disjoint edges such that $\sign(ij,k\ell)\ne 0$, with $i<j$ and $k<\ell$. Finally, the evaluation map $E:\ss\to\real^{\xx}$ sends $P$ to the vector
$$
E(P)=\big(P(i,j,k,\ell)\big)_{(ij,k\ell)\in\xx}.
$$


\subsection{The strategy and the key objects}\label{sub:strat}

In a nutshell the proof has two main steps. First we show that every two vectors of $E(\ss)$ are orthogonal with respect to a weighted version of the ordinary dot product on $\real^{\xx}$, in which some of the weights are negative. For $a,b\in\real^{\xx}$ we let
\begin{equation}\label{eq:iform}
\langle a,b\rangle_\mu\mathrel{:=}\sum_{(ij,k\ell)\in\xx}\mu(ij,k\ell)\,a_{(ij,k\ell)}\,b_{(ij,k\ell)},
\end{equation}
where the coefficients $\mu(ij,k\ell)$ are nonzero real numbers that we define below. We say that two vectors $a$ and $b$ in $\real^{\xx}$ are {\em orthogonal} if $\langle a,b\rangle_\mu=0$. 

As we mentioned above, in the first step we show that every two vectors of $E(\ss)$ are orthogonal with respect to $\langle\cdot,\cdot\rangle_\mu$. The second step consists of proving that every subspace with this property has dimension at most $\frac{1}{2}|\xx|$. This last fact is totally elementary, and uses only that the coefficients $\mu(ij,k\ell)$ are nonzero.

The first step is then the more difficult one, and it entirely relies on the specific coefficients $\mu(ij,k\ell)$. With an eye on the Main Identity we let
\begin{equation}\label{eq:weightc}
\mu(ij,k\ell)\mathrel{:=}\sign(ij,k\ell)\,\omega_i\omega_j\omega_k\omega_\ell\,L(i,j,k,\ell),
\end{equation}
where 
\begin{equation}\label{eq:Delta}
L(x_1,x_2,y_1,y_2)\mathrel{:=}(x_1-x_2)(y_1-y_2)(x_1-y_1)(x_1-y_2)(x_2-y_1)(x_2-y_2).
\end{equation}

As we are about to see, this choice of $\mu(ij,k\ell)$ accomplishes the claim that every two vectors of $E(\ss)$ are orthogonal with respect to $\langle\cdot,\cdot\rangle_\mu$.

We note that $\langle\cdot,\cdot\rangle_\mu$ is not a dot product in the usual sense, as the definition of sign plus the identity $L(k,\ell,i,j)=L(i,j,k,\ell)$ imply that $\mu(k\ell,ij)=-\mu(ij,k\ell)$. Thus half of the coefficients $\mu(ij,k\ell)$ are negative, and a nonzero vector may be orthogonal to itself (which is precisely what happens with the vectors of $E(\ss)$).

\subsection{Proof of Proposition~\ref{pro:iso}}\label{sub:isoproof}

In the proof we use the following elementary linear algebra fact.

\begin{lemma}\label{lem:halfdim}
Let $c_1,\ldots,c_m$ be nonzero real numbers, and let $W$ be a subspace of $\real^m$ such that
$$
c_1a_1b_1+c_2a_2b_2+\cdots+c_ma_mb_m=0\quad\text{for all $a,b\in W$.}
$$
Then $\dim W\le\tfrac{1}{2}m$.
\end{lemma}

\begin{proof}
Let $W'$ be the set of the vectors $(c_1b_1,\ldots,c_mb_m)$ with $b\in W$. Since the numbers $c_t$ are nonzero, $W'$ is a subspace of $\real^m$ of the same dimension as $W$. By hypothesis, every vector of $W'$ has ordinary dot product zero with every vector of $W$, and so $W'$ is contained in the orthogonal complement of $W$, which has dimension $m-\dim W$. Therefore $\dim W\le m-\dim W$, that is, $\dim W\le\tfrac{1}{2}m$.
\end{proof}

\begin{proof}[Proof of Proposition~\ref{pro:iso}]
As we mentioned above, we start by showing that ($*$) {\em for all $P,P'\in\ss$, the vectors $E(P)$ and $E(P')$ are orthogonal.}

To show ($*$) let $P,P'$ be polynomials in $\ss$, and let $A=P\cdot P'\cdot L$, where $L$ is defined in~\eqref{eq:Delta}. We start by noting that $A$ is an antisymmetric pair polynomial, so we can apply the Main Identity to it. To see this we note that since $P$ and $P'$ have degree at most $r-2$ in each variable, and $L$ has degree $3$ in each variable, it follows that $A$ has degree at most $2(r-2)+3=2r-1=n-2$ in each variable, and so $A$ is a pair polynomial. We finally observe that $A$ is antisymmetric because $P$ and $P'$ are symmetric and $L$ is antisymmetric.

Thus $A$ is indeed an antisymmetric pair polynomial, and so we can apply the Main Identity (Theorem~\ref{thm:master}) to it. We first note that the rotations side is $0$. Indeed, each of the terms of this side involves a value $A(i,j,i,k)$, that is, a term in which the first and the third variables take the same value $i$, and since $L$ has the factor $x_1-y_1$ it follows that $L(i,j,i,k)=0$, and so $A(i,j,i,k)=0$.

The Main Identity then implies that the quads side is also zero. This side is
$$
\sum_{(ij,k\ell)\in\xx}\sign(ij,k\ell)\,\omega_i\omega_j\omega_k\omega_\ell\,L(i,j,k,\ell)\,P(i,j,k,\ell)\,P'(i,j,k,\ell),
$$
and the definitions of $\mu(ij,k\ell)$ in~\eqref{eq:weightc} and of $\langle\cdot,\cdot\rangle_\mu$ in~\eqref{eq:iform} imply that this is exactly $\langle E(P),E(P')\rangle_\mu$. Thus ($*$) holds.

In the second and final step of the proof we apply Lemma~\ref{lem:halfdim} to the subspace $W=E(\ss)$ of $\real^{\xx}$, with $m=|\xx|$ and with the coefficients $\mu(ij,k\ell)$ as the numbers $c_1,\ldots,c_m$ (we list the elements of $\xx$ in any order). We note that all these coefficients are nonzero, as required in the lemma. Indeed, $\sign(ij,k\ell)\ne 0$ for $(ij,k\ell)\in\xx$, the weights $\omega_i$ are all nonzero, and $L(i,j,k,\ell)\ne 0$ since $i,j,k,\ell$ are distinct. 

By the definition of $\langle\cdot,\cdot\rangle_\mu$, the hypothesis of Lemma~\ref{lem:halfdim} for $W=E(\ss)$ is exactly ($*$). Therefore Lemma~\ref{lem:halfdim} applies, and so $\dim E(\ss)\le\tfrac{1}{2}|\xx|$.
\end{proof}


\section{Proof of Lemma~\ref{lem:detector}}\label{sec:interp}

In the context of Lemma~\ref{lem:detector}, we recall that $n=2r+1$ and that $\omega_j=1/\prod_{k\in[n],\,k\ne j}(j-k)$ for every $j\in[n]$. For easy reference let us highlight the main objects in the lemma.

\begin{defbox}{Setting for this section}
Throughout this section:
\begin{itemize}
\item $i\in[n]$ is fixed.
\item $N:=[n]\setminus\{i\}$.
\item $\rho$ is the sign function of a fixed linear order of the $2r$ integers in $N$.
\end{itemize}
\end{defbox}

\subsection{Roadmap of the proof of Lemma~\ref{lem:detector}}\label{sub:roadmap10}

In Section~\ref{sub:rewriting} we reformulate Lemma~\ref{lem:detector} in terms of a weighted sum $\llbracket F,B\rrbracket$ that resembles a dot product, and with this analogy in mind we say that $F$ and $B$ are orthogonal if $\llbracket F,B\rrbracket=0$. In this reformulation we also use an extension of the weights $\omega_i$ to a more general version that allows arbitrary finite sets of real numbers.

Section~\ref{sub:twofacts} contains the workhorse behind the proof of Lemma~\ref{lem:detector}, a technical lemma (namely Lemma~\ref{lem:twotri}) that says that if a polynomial in two variables of low degree vanishes at a certain set $\mathcal A_r$ of $r(r-1)$ points then it must be the zero polynomial.

Finally, in Section~\ref{sub:proofdetector} we prove Lemma~\ref{lem:detector}. The key idea is to choose $B$ cleverly, so that $\llbracket F,B\rrbracket$ becomes a sum of squares. Since $\llbracket F,B\rrbracket=0$ it will follow that each of these squares is zero, and therefore $F$ vanishes on $\mathcal A_r$. Using Lemma~\ref{lem:twotri} it will then follow that $F$ is the zero polynomial.

\subsection{Restating Lemma~\ref{lem:detector} in terms of a weighted sum}\label{sub:rewriting}

We start by extending the weights $\omega_1,\ldots,\omega_n$ of Section~\ref{sec:conv} to arbitrary finite sets of numbers. Let $Y$ be a finite set of real numbers. For each $j\in Y$ we let
\begin{equation}\label{eq:omegaNdef}
\omega_Y(j)\mathrel{:=}\frac{1}{\prod_{k\in Y,\,k\ne j}(j-k)}.
\end{equation}
Thus for $Y=[n]$ these are the weights $\omega_j$ defined in Section~\ref{sec:conv} that we have used up to this point, that is, $\omega_{[n]}(j)=\omega_j$.

We now note that Lemma~\ref{lem:lagrange} holds for these more general weights as well.

\begin{lemma}\label{lem:lagrange-gen}
Let $Y$ be a finite set of real numbers. Then for every polynomial $G$ of degree at most $|Y|-2$,
$$
\sum_{j\in Y}\omega_{Y}(j)\,G(j)=0.
$$
\end{lemma}

\begin{proof}[Proof of Lemma~\ref{lem:lagrange-gen}]
The proof is virtually identical to the proof of Lemma~\ref{lem:lagrange}. It suffices to replace $[n]$ with $Y$, $n$ with $|Y|$, and $\omega_i$ with $\omega_{Y}(i)$.
\end{proof}

In this section we mostly use the case $Y=N$, in which $\omega_N(j)$ is similar to $\omega_j$, but computed as if $i$ were not in $[n]$. That is,
\begin{equation}\label{eq:omegaN}
\omega_N(j)=\omega_j\,(j-i)\quad\text{for every $j\in N$.}
\end{equation}
For polynomials $F$ and $B$ in two variables, we define the weighted sum $\llbracket F,B\rrbracket$ by
\begin{equation}\label{eq:pairing}
\llbracket F,B\rrbracket\mathrel{:=}\sum_{j\in N}\omega_N(j)\sum_{k\in N}\omega_N(k)\,\tfrac{1}{2}\rho(j,k)\,(j-k)\,F(j,k)\,B(j,k).
\end{equation}

Since $\llbracket F,B\rrbracket$ resembles a dot product, it will be useful to borrow the following term from linear algebra.

\begin{defbox}{Orthogonality}
Throughout this section we say that two polynomials $F$ and $B$ in two variables are {\em orthogonal} if $\llbracket F,B\rrbracket=0$.
\end{defbox}

In view of~\eqref{eq:omegaN}, the product $\omega_j\omega_k\,(j-i)(k-i)$ in~\eqref{eq:detector} is $\omega_N(j)\,\omega_N(k)$, and so the sum in~\eqref{eq:detector} is $2\llbracket F,B\rrbracket$. Therefore Lemma~\ref{lem:detector} can be restated as follows.

\theoremstyle{plain}
\newtheorem*{restateddetector}{Lemma~\ref{lem:detector}}
\begin{restateddetector}[restated in terms of $\llbracket\cdot,\cdot\rrbracket$]
Let $F(x_2,y_2)$ be a polynomial of degree at most $r-2$ in each variable. If $F$ is orthogonal to every polynomial $B(x_2,y_2)$ of degree at most $r-1$ in each variable, then $F$ is the zero polynomial.
\end{restateddetector}

\subsection{A vanishing lemma for polynomials in two variables}\label{sub:twofacts}

Lemma~\ref{lem:twotri} below is a two-variable version of the elementary fact that a polynomial of degree at most $d$ in one variable is determined by its values at $d+1$ distinct points. It says that a polynomial in two variables of low degree is determined by its values at a particular set of $r(r-1)$ points. 

We use $\pp_{d_1,d_2}$ to denote the space of real polynomials $F(x,y)$ in two variables $x$ and $y$, of degree at most $d_1$ in $x$ and at most $d_2$ in $y$. (For the polynomials $F(x_2,y_2)$ and $B(x_2,y_2)$ of Lemma~\ref{lem:detector}, the variables $x$ and $y$ stand for $x_2$ and $y_2$.) For instance,
$$
F(x,y)=3x^2y-xy+5
$$
belongs to $\pp_{2,1}$. We note that $\dim\pp_{d_1,d_2}=(d_1+1)(d_2+1)$.

\begin{figure}[ht!]
\centering
\includegraphics{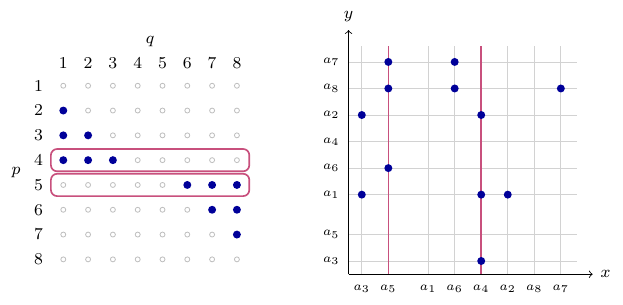}
\caption{Left: the pattern $\mathcal L_4$ for eight nodes. The index pair $(p,q)$ is drawn in row $p$ and column $q$, and the pairs in $\mathcal L_4$ are the filled dots. The rows $p=4$ and $p=5$, framed, have $r-1=3$ pairs each. They are the rows used in the proof of Lemma~\ref{lem:twotri}. Right: the corresponding set $\mathcal A_4$ of points of the plane, for the list $a_1,\ldots,a_8=3,6,\tfrac{1}{2},5,\tfrac{3}{2},4,8,7$. Each point $(a_p,a_q)$ lies on the vertical line $x=a_p$ and on the horizontal line $y=a_q$. The framed rows $p=4$ and $p=5$ on the left correspond to the vertical lines $x=a_4$ and $x=a_5$ on the right (thick), each of which contains $r-1=3$ points of $\mathcal A_4$.}
\label{fig:twotri}
\end{figure}

Let $a_1,\ldots,a_{2r}$ be distinct real numbers, which we call {\em nodes}, listed in any order. We are interested in the set of index pairs
$$
\mathcal L_r=\{(p,q):1\le q<p\le r\}\cup\{(p,q):r<p<q\le 2r\}.
$$
We refer the reader to the left-hand side of Figure~\ref{fig:twotri}, where we illustrate $\mathcal L_4$. Loosely speaking, $\mathcal L_r$ consists of two ``grid triangles''. Each of them has $\binom{r}{2}$ elements, and so $|\mathcal L_r|=2\binom{r}{2}=r(r-1)$.

The corresponding set of points of the plane is
$$
\mathcal A_r=\{(a_p,a_q):(p,q)\in\mathcal L_r\}.
$$
Since the nodes are distinct these points are distinct, and so $|\mathcal A_r|=|\mathcal L_r|=r(r-1)$. The right-hand side of Figure~\ref{fig:twotri} shows $\mathcal A_4$ for the choices of the nodes $a_1,\ldots,a_8$ listed in the caption.

We say that a polynomial $F$ in two variables {\em vanishes on $\mathcal A_r$} if $F(a_p,a_q)=0$ for every $(p,q)\in\mathcal L_r$.

\begin{lemma}\label{lem:twotri}
If a polynomial in $\pp_{r-1,r-2}$ or in $\pp_{r-2,r-1}$ vanishes on $\mathcal A_r$, then it is the zero polynomial.
\end{lemma}

\begin{proof}[Proof of Lemma~\ref{lem:twotri}]
We actually prove the following stronger statement: the lemma holds for every integer $r\ge 2$, and not only for $r=(n-1)/2$.

We proceed by induction on $r$, proving the statement for both spaces $\pp_{r-1,r-2}$ and $\pp_{r-2,r-1}$ at once. We start with the base case $r=2$, for the spaces $\pp_{1,0}$ and $\pp_{0,1}$. In this case $\mathcal A_2$ consists of the two points $(a_2,a_1)$ and $(a_3,a_4)$, whose first coordinates are distinct, and whose second coordinates are also distinct. Now a polynomial in $\pp_{1,0}$ is of the form $c\,x+d$ (it does not involve the second variable), and if it vanishes at two points with distinct first coordinates then it is the zero polynomial. An analogous remark holds for polynomials in $\pp_{0,1}$. This takes care of the case $r=2$.

For the inductive step we let $r>2$ be an integer, assume that the statement holds for $r-1$, that is, for $\pp_{r-2,r-3}$ and $\pp_{r-3,r-2}$, and prove that then it holds for $\pp_{r-1,r-2}$ and $\pp_{r-2,r-1}$. Thus we let $F$ be a polynomial in $\pp_{r-1,r-2}$ or in $\pp_{r-2,r-1}$ that vanishes on $\mathcal A_r$, and our goal is to show that $F$ is the zero polynomial. We deal with the two spaces separately.

Suppose first that $F\in\pp_{r-1,r-2}$. We note that each of the rows $p=r$ and $p=r+1$ contains $r-1$ pairs of $\mathcal L_r$, as we illustrate in Figure~\ref{fig:twotri}. The row $p=r$ contains the pairs $(r,q)$ with $1\le q\le r-1$, and the row $p=r+1$ contains the pairs $(r+1,q)$ with $r+2\le q\le 2r$.

Consider the one-variable polynomial $F(a_r,y)$. Its degree is at most $r-2$, and it vanishes at the $r-1$ distinct points $y=a_q$ with $1\le q\le r-1$. Thus $F(a_r,y)$ is the zero polynomial, since a nonzero polynomial in one variable of degree at most $r-2$ has at most $r-2$ distinct roots. In a totally analogous manner we conclude that $F(a_{r+1},y)$ is the zero polynomial. Therefore $x-a_r$ and $x-a_{r+1}$ divide $F(x,y)$ (indeed, if we regard $F(x,y)$ as a polynomial in $y$ whose coefficients are polynomials in $x$, then each of these coefficients vanishes at $x=a_r$ and at $x=a_{r+1}$), and so there exists a polynomial $G\in\pp_{r-3,r-2}$ such that
$$
F(x,y)=(x-a_r)(x-a_{r+1})\,G(x,y).
$$

We now remove the two nodes $a_r$ and $a_{r+1}$ from the list $a_1,\ldots,a_{2r}$, and let $a'_1,\ldots,a'_{2r-2}$ be the remaining $2r-2$ nodes, listed in the same order. A glance at Figure~\ref{fig:twotri} reveals that the pairs of $\mathcal L_r$ that survive the removal of the rows and columns $r$ and $r+1$ are exactly the pairs of $\mathcal L_{r-1}$ for this new list. Moreover, the factor $(x-a_r)(x-a_{r+1})$ is nonzero at each of these surviving points, since their first coordinates are neither $a_r$ nor $a_{r+1}$. 

Since $F$ vanishes on $\mathcal A_r$ then in particular it vanishes at these surviving points, and since the factor $(x-a_r)(x-a_{r+1})$ is nonzero at each of them it follows that $G$ vanishes on the set $\mathcal A_{r-1}$ of points that corresponds to this new list. Since $G\in\pp_{r-3,r-2}=\pp_{(r-1)-2,(r-1)-1}$ the induction hypothesis (for the second space) implies that $G$ is the zero polynomial. Therefore $F$ is the zero polynomial, as required.

Suppose finally that $F\in\pp_{r-2,r-1}$. The argument is totally analogous, using columns instead of rows. Each of the columns $q=1$ and $q=2r$ contains exactly $r-1$ pairs of $\mathcal L_r$, namely the pairs $(p,1)$ with $2\le p\le r$, and the pairs $(p,2r)$ with $r+1\le p\le 2r-1$.

Thus the one-variable polynomials $F(x,a_1)$ and $F(x,a_{2r})$, each of degree at most $r-2$, vanish at $r-1$ distinct points, and so they are both the zero polynomial. Hence $F(x,y)=(y-a_1)(y-a_{2r})\,G(x,y)$ for some $G\in\pp_{r-2,r-3}$.

Removing the nodes $a_1$ and $a_{2r}$, the surviving pairs of $\mathcal L_r$ turn out to be $\mathcal L_{r-1}$ for the remaining nodes, and $G$ vanishes on the corresponding set of points $\mathcal A_{r-1}$. Since $G\in\pp_{(r-1)-1,(r-1)-2}$, the induction hypothesis (for the first space) implies that $G$ is the zero polynomial, and so $F$ is the zero polynomial.
\end{proof}


\subsection{Proof of Lemma~\ref{lem:detector}}\label{sub:proofdetector}
We start by noting that in order to prove the lemma it suffices to show the following.

\vglue 0.3 cm
\noindent($\ast$)\,\,(implies Lemma~\ref{lem:detector}) {\em If $F\in\pp_{r-1,r-2}$ is orthogonal to every polynomial in $\pp_{r-2,r-1}$, then $F$ is the zero polynomial.}
\vglue 0.3 cm

To see that ($\ast$) indeed implies Lemma~\ref{lem:detector}, we recall from Section~\ref{sub:rewriting} that Lemma~\ref{lem:detector} is equivalent to the following statement: {\em if $F$ has degree at most $r-2$ in each variable, and it is orthogonal to every polynomial of degree at most $r-1$ in each variable, then $F$ is the zero polynomial}. In the notation $\pp_{d_1,d_2}$ of Section~\ref{sub:twofacts}, the hypothesis is that $F\in\pp_{r-2,r-2}$ is orthogonal to every polynomial in $\pp_{r-1,r-1}$.

Suppose then that $F\in\pp_{r-2,r-2}$ is orthogonal to every polynomial in $\pp_{r-1,r-1}$. Note that in particular $F\in\pp_{r-1,r-2}$, since $\pp_{r-2,r-2}\subseteq\pp_{r-1,r-2}$. Moreover $F$ is orthogonal to every polynomial in $\pp_{r-2,r-1}$, since $\pp_{r-2,r-1}\subseteq\pp_{r-1,r-1}$. Thus $F$ satisfies the hypothesis of ($\ast$), and under the assumption that ($\ast$) holds it follows that $F$ is the zero polynomial, as required.

We now prove ($\ast$). For the rest of the proof $F$ is a polynomial in $\pp_{r-1,r-2}$ and $B$ is a polynomial in $\pp_{r-2,r-1}$. 

We recall that the orthogonality between $F$ and $B$ involves the weight function $\omega_N$ as well as $\rho$. We recall that $\rho$ is the sign function of a fixed linear order of the $2r$ integers in $N$. For the rest of the proof, we let this fixed linear order be $j_1,j_2,\ldots,j_{2r}$. Thus $\rho(j_p,j_q)=+1$ if $p>q$, $\rho(j_p,j_q)=-1$ if $p<q$, and $\rho(j_p,j_p)=0$.

In a nutshell, the strategy to prove ($\ast$) is to show that $F(j_p,j_q)=0$ for every $(p,q)\in\mathcal L_r$, that is, that $F$ vanishes on $\mathcal A_r$ for the list $j_1,\ldots,j_{2r}$. Once we show this we will be done, as then Lemma~\ref{lem:twotri} implies that $F$ is the zero polynomial.

To help comprehension we break the proof of ($\ast$) into three steps.

\medskip
\noindent{\sc Step 1}.\enspace{\em The value of $\llbracket F,B\rrbracket$ does not change if in~\eqref{eq:pairing} we replace $\tfrac{1}{2}\rho(j,k)$ by $\tfrac{1}{2}\rho(j,k)+\lambda(j)$, where $\lambda$ is an arbitrary real function on $N$.}\par\medskip

It is straightforward to check that the difference in $\llbracket F,B\rrbracket$ if we replace $\tfrac{1}{2}\rho(j,k)$ by $\tfrac{1}{2}\rho(j,k)+\lambda(j)$ is
$$
\sum_{j\in N}\lambda(j)\,\omega_N(j)\sum_{k\in N}\omega_N(k)\,(j-k)\,F(j,k)\,B(j,k).
$$
For each fixed $j\in N$, the expression $(j-k)F(j,k)B(j,k)$ is a polynomial in $k$ of degree at most $1+(r-2)+(r-1)=2r-2=|N|-2$, and so by Lemma~\ref{lem:lagrange-gen} the summation over $k\in N$ is zero. Thus Step 1 follows.

\medskip
\noindent{\sc Step 2}.\enspace{\em The weighted sum $\llbracket F,B\rrbracket$ can be written as a sum over $\mathcal L_r$. More precisely,
\begin{equation}\label{eq:overL}
\llbracket F,B\rrbracket=\sum_{(p,q)\in\mathcal L_r}\eta(p,q)\,(j_p-j_q)\,\omega_N(j_p)\,\omega_N(j_q)\,F(j_p,j_q)\,B(j_p,j_q),
\end{equation}
where $\eta(p,q)=+1$ if $(p,q)$ is in the first triangle of $\mathcal L_r$, and $\eta(p,q)=-1$ if it is in the second triangle.}\par\medskip

We start by writing $\llbracket F,B\rrbracket$ as a sum over all pairs $(p,q)$. Since $N=\{j_1,\ldots,j_{2r}\}$, the definition of $\llbracket F,B\rrbracket$ (see~\eqref{eq:pairing}) reads
\begin{equation}\label{eq:allpq}
\llbracket F,B\rrbracket=\sum_{p=1}^{2r}\omega_N(j_p)\sum_{q=1}^{2r}\omega_N(j_q)\,\tfrac{1}{2}\rho(j_p,j_q)\,(j_p-j_q)\,F(j_p,j_q)\,B(j_p,j_q).
\end{equation}
This sum looks like the sum in~\eqref{eq:overL}, but there are two differences. First, this sum is over all the pairs $(p,q)$, not only over the pairs in $\mathcal L_r$. Second, the coefficients are $\tfrac{1}{2}\rho(j_p,j_q)$ instead of $\pm 1$. We remove both differences with a suitable application of Step~1.

For simplicity we define $\mathbf 1[p>r]$ by saying that it is $1$ if $p>r$ and $0$ otherwise.

Using Step~1 with $\lambda(j_p)=\tfrac{1}{2}-\mathbf 1[p>r]$ the sum in~\eqref{eq:allpq} does not change if we replace each coefficient $\tfrac{1}{2}\rho(j_p,j_q)$ by
$$
\eta(p,q)\mathrel{:=}\tfrac{1}{2}\rho(j_p,j_q)+\tfrac{1}{2}-\mathbf 1[p>r].
$$
That is,
\begin{equation}\label{eq:allpqeps}
\llbracket F,B\rrbracket=\sum_{p=1}^{2r}\omega_N(j_p)\sum_{q=1}^{2r}\omega_N(j_q)\,\eta(p,q)\,(j_p-j_q)\,F(j_p,j_q)\,B(j_p,j_q).
\end{equation}

The terms of~\eqref{eq:allpqeps} with $p=q$ are zero, because they have the factor $j_p-j_q$. For $p\ne q$ we have $\rho(j_p,j_q)=+1$ if $p>q$ and $\rho(j_p,j_q)=-1$ if $p<q$, and so $\eta(p,q)=\tfrac{1}{2}\rho(j_p,j_q)+\tfrac{1}{2}$ is $1$ or $0$ if $p\le r$, while $\eta(p,q)=\tfrac{1}{2}\rho(j_p,j_q)-\tfrac{1}{2}$ is $0$ or $-1$ if $p>r$. In other words, for $p\ne q$ we have
$$
\eta(p,q)=
\begin{cases}
+1 & \text{if $q<p\le r$, that is, if $(p,q)$ is in the first triangle of $\mathcal L_r$,}\\
-1 & \text{if $r<p<q$, that is, if $(p,q)$ is in the second triangle of $\mathcal L_r$,}\\
\phantom{+}0 & \text{otherwise.}
\end{cases}
$$

Therefore the only terms of~\eqref{eq:allpqeps} that survive are those with $(p,q)\in\mathcal L_r$, with $\eta(p,q)=+1$ on the first triangle and $\eta(p,q)=-1$ on the second one. This is exactly~\eqref{eq:overL}, and so Step~2 follows.

\medskip
\noindent{\sc Step 3}.\enspace{\em If $F$ is orthogonal to every $B\in\pp_{r-2,r-1}$, then $F$ is the zero polynomial. That is, ($\ast$) holds.}\par\medskip

Suppose that $F$ is orthogonal to every $B\in\pp_{r-2,r-1}$. For brevity, for each $(p,q)\in\mathcal L_r$ we let $c_{pq}=\eta(p,q)\,(j_p-j_q)\,\omega_N(j_p)\,\omega_N(j_q)$, so that~\eqref{eq:overL} reads
\begin{equation}\label{eq:overLc}
\llbracket F,B\rrbracket=\sum_{(p,q)\in\mathcal L_r}c_{pq}\,F(j_p,j_q)\,B(j_p,j_q).
\end{equation}
We emphasize that $c_{pq}\ne 0$ for every $(p,q)\in\mathcal L_r$, since $\eta(p,q)=\pm 1$, $j_p\ne j_q$, and all the weights $\omega_N(j_p)$ are nonzero.

We now apply Lemma~\ref{lem:twotri} to $j_1,\ldots,j_{2r}$ (that is, with $a_p=j_p$ for $p=1,\ldots,2r$). Evaluation at the points of $\mathcal A_r$ is a linear map from $\pp_{r-2,r-1}$ to $\real^{\mathcal A_r}$. By Lemma~\ref{lem:twotri} this map is injective, since the zero polynomial is the only polynomial mapped to the zero vector. Since both spaces have dimension $r(r-1)=|\mathcal A_r|$, this map is also surjective. Thus we can find a polynomial $B\in\pp_{r-2,r-1}$ that takes any values we want at the points of $\mathcal A_r$.

With this freedom we choose $B$ so that $B(j_p,j_q)=c_{pq}F(j_p,j_q)$ for every $(p,q)\in\mathcal L_r$. With this choice each term in the sum in~\eqref{eq:overLc} is $(c_{pq}F(j_p,j_q))^2$, and so $\llbracket F,B\rrbracket$ is a sum of squares. Since $F$ is orthogonal to $B$, this sum of squares is zero, and so each of the squared numbers is zero. Since $c_{pq}\ne 0$, it follows that $F(j_p,j_q)=0$ for every $(p,q)\in\mathcal L_r$, that is, $F$ vanishes on $\mathcal A_r$. Since $F\in\pp_{r-1,r-2}$, Lemma~\ref{lem:twotri} implies that $F$ is the zero polynomial.\qed


\section{Concluding remarks}\label{sec:concluding}

As we mentioned in Section~\ref{sub:ourmain}, Dan Archdeacon proposed (in our terminology) the following conjecture in his list of problems in topological graph theory~\cite{Archdeacon99} (see also~\cite{AMRS}, where his conjecture is reproduced). 

\begin{conjecture}[Archdeacon]\label{con:archdeacon}
If $\Pi$ is any rotation system of $K_n$ then
$$
\crq{\Pi}+\nr(\Pi)\ge Z(n).
$$
\end{conjecture}

Evidently Conjecture~\ref{con:archdeacon} is stronger than our Theorem~\ref{thm:odd}. Chikkatur and Zeng~\cite{ChikkaturZeng} verified the conjecture for $n\le 10$ by computer, and proved that $\crq{\Pi}+\nr(\Pi)\ge (8/9-o(1))Z(n)$.

We finally point out that Balko~\cite{balkohill} and Schaefer and \v{S}tefankovi\v{c}~\cite{SchaeferStefankovic} have also written expositions of the proof in~\cite{openai}, both of which contain further mathematical results. Let us briefly discuss the purpose and scope of these works.

As we were writing this paper it came to our attention that Martin Balko was also working on an exposition of OpenAI's proof of the Harary--Hill conjecture. We got in touch, and before taking a look at each other's manuscripts we briefly considered the possibility of writing a joint paper. However, when we realized that our aims were quite different we decided that the best course of action would be to continue to work separately, keeping each other up to date in this process. 

It is worth mentioning that in the last section of his paper Balko~\cite{balkohill} shows that the Harary--Hill conjecture combined with known reductions easily determines the crossing numbers of several other families of graphs that have been considered in the literature. These include the complete graph $K_{2k}$ with a perfect matching or a single edge removed, the Cartesian products of complete graphs with paths and with even cycles, and two extensions of Sierpi\'nski graphs. Balko also shows that $\ocr{K_n}$ and $\pcr{K_n}$ (the odd crossing number and the pair crossing number of $K_n$, when restricted to drawings in which adjacent edges do not cross) are both equal to $Z(n)$. 

Very recently Schaefer and \v{S}tefankovi\v{c}~\cite{SchaeferStefankovic} uploaded another exposition of OpenAI's proof. Their nicely streamlined proof is considerably strengthened by showing in addition that the independent odd crossing number $\iocr{K_n}$ of $K_n$ is $Z(n)$. Since several other crossing number variants, including the odd crossing number, the algebraic crossing number, and the pair crossing number, lie between the independent odd crossing number and the crossing number (see~\cite{SchaeferSurvey} for these facts, and for the definitions of these variants), their result determines these crossing number variants for complete graphs as well.

Remarkably, their method also establishes that the independent odd crossing number of $K_{m,n}$ is $Z(m,n)$. This determines the corresponding crossing number variants for complete bipartite graphs and, in particular, establishes Zarankiewicz's conjecture, also recently proved by OpenAI~\cite{openai-bip}.


\section{Acknowledgments}

We thank Martin Balko and Ruy Fabila for useful discussions. We also thank Alan Arroyo, Dan McQuillan, Bruce Richter, and Matthew Sullivan for engaging discussions on the rotation systems of complete graphs. In particular, the insights on how the orientations of the vertices in a quad lead to the signs of a quad can be traced back to Sullivan's lucid explanations during our group work sessions.


\bibliographystyle{amsrefs}
\bibliography{refs}


\end{document}